\pdfoutput=1
\documentclass[reqno]{amsart}

\usepackage{caption}
\usepackage{verbatim}
\usepackage{mathrsfs}
\usepackage{stmaryrd}

\usepackage[paperwidth=7in,paperheight=10in,text={5in,8in},centering]{geometry}

\usepackage{latexsym}
\usepackage{amsfonts}
\usepackage{amsmath}
\usepackage{amssymb}
\usepackage{graphicx}
\usepackage{comment}
\usepackage{soul}
\usepackage{color}

\usepackage[numbers,square,comma]{natbib}
\usepackage{hyperref}
\usepackage{titlesec}

\hypersetup{
    colorlinks=true,
    linkcolor=blue,
    citecolor=blue,
    pdfborder={0 0 0}
}

\newcounter{saveeqn}

\newcommand{\oo}{\infty}
\newcommand{\pp}[2]{\frac{\partial #1}{\partial #2}}

\titleformat{\section}[runin]
{\normalfont\bfseries}{\thesection .}{1ex}{}[]
  
\titleformat{\subsection}[runin]
{\normalfont\itshape}{\thesubsection}{1ex}{}

\titleformat{\subsubsection}[runin]
{\normalfont\itshape}{\thesubsubsection}{1ex}{}

\theoremstyle{definition}

\newcommand{\sD}{\mathscr{D}}
\newcommand{\sN}{\mathscr{N}}

\newcommand{\e}[1]{{(#1)}}
\newcommand{\veps}{\varepsilon}
\newcommand{\jd}{\displaystyle}
\newcommand{\jt}{\textstyle}
\newcommand{\opn}[1]{\operatorname{#1}}
\newcommand{\der}[2]{\frac{\partial{#1}}{\partial{#2}}}

\title[Computation of Dirichlet--Neumann Operators]{Comparison of Five
  Methods of Computing \\ the Dirichlet--Neumann Operator for \\
  the Water Wave Problem}

\author[J. Wilkening]{Jon Wilkening}
\address{Department of Mathematics, University of California, Berkeley CA 94720-3840}
\author[V. Vasan]{Vishal Vasan}
\address{Department of Mathematics, Pennsylvania State University, University Park PA 16802}

\date{\today}

\begin{document}

\begin{abstract}
  We compare the effectiveness of solving Dirichlet--Neumann problems
  via the Craig--Sulem (CS) expansion, the Ablowitz-Fokas-Musslimani
  (AFM) implicit formulation, the dual AFM formulation (AFM$^*$), a
  boundary integral collocation method (BIM), and the transformed
  field expansion (TFE) method.  The first three methods involve
  highly ill-conditioned intermediate calculations that we show can be
  overcome using multiple-precision arithmetic.  The latter two
  methods avoid catastrophic cancellation of digits in intermediate
  results, and are much better suited to numerical computation.

  For the Craig--Sulem expansion, we explore the cancellation of terms
  at each order (up to 150th) for three types of wave profiles, namely
  band-limited, real-analytic, or smooth.  For the AFM and AFM$^*$
  methods, we present an example in which representing the Dirichlet
  or Neumann data as a series using the AFM basis functions is
  impossible, causing the methods to fail. The example involves
  band-limited wave profiles of arbitrarily small amplitude, with
  analytic Dirichlet data.  We then show how to regularize the AFM and
  AFM$^*$ methods by over-sampling the basis functions and using the
  singular value decomposition or QR-factorization to orthogonalize
  them.  Two additional examples are used to compare all five methods
  in the context of water waves, namely a large-amplitude standing
  wave in deep water, and a pair of interacting traveling
  waves in finite depth.
\end{abstract}

\maketitle

\section{Introduction}

The water wave equations, also known as Euler's equations for
inviscid, irrotational waves, describe the motion of the free surface
of an ideal fluid. For periodic waves in two dimensions over a flat
surface in the absence of surface tension, they are given as
follows \cite{Kundu}:
\begin{alignat*}{2}
    \phi_{xx} + \phi_{yy}  &=0,\quad &-h&<y<\eta,\\
    \eta_{t} + \phi_{x}\eta_{x} &= \phi_{y},\quad &y&=\eta(x,t),\\
    \phi_{t} + \frac 1 2 \phi_{x}^{2} + \frac 1 2 \phi_{y}^{2} + g\eta &= 0,
    \quad &y&=\eta(x,t), \\
    \phi_{y} &= 0,\quad &y&=-h.
\end{alignat*}
Here $\phi$ is the velocity potential (related to the
fluid velocity by $\mathbf{u}=\nabla\phi$), $\eta$ is the graph of the 
free surface of the water, $g$ is the acceleration of gravity,
$h$ is the depth of the undisturbed fluid, and subscripts denote
partial derivatives. Since we are interested in 
periodic waves, we consider periodic boundary conditions in the 
horizontal direction. We note that the above equations represent 
a free-boundary value problem for Laplace's equation, \emph{i.e.}~both 
$\eta$ and $\phi$ are unknowns. These equations readily generalize 
to three-dimensional fluids (with two-dimensional surfaces);
however, in the present work, we limit ourselves to one-dimensional 
surfaces.

The water-wave equations may be reformulated in terms of only the 
surface variables $\eta(x,t)$ and $q=\phi(x,\eta,t)$ as \cite{Zakharov,CraigSulem}
\begin{align*}
    \frac{\partial \eta}{\partial t} &= G(\eta)q,\\
    \frac{\partial q}{\partial t} &= -g\eta - \frac{q_{x}^{2}}{2} +
    \frac{(G(\eta)q+\eta_{x}q_{x})^{2}}{2(1+\eta_{x}^{2})},
\end{align*}
where $G(\eta)$ represents the Dirichlet-Neumann operator (DNO) defined 
as 
\[G(\eta)q = \phi_{y}(x,\eta) - \eta_{x}\phi_{x}(x,\eta),\] where $\phi${}
 is the solution to the following boundary-value problem
\begin{alignat*}{2}
    \phi_{xx} + \phi_{yy}  &=0, &  -h&<y<\eta,\\
    \phi(x,y) &= q(x), \quad & y &= \eta, \\
    \phi_{y} &= 0, &  y&=-h.
\end{alignat*}
Thus $G(\eta)$ maps the given Dirichlet data $q$ to the associated 
Neumann data at the free surface. Consequently, to evolve the surface 
variables in time using some numerical scheme, we require the solution 
to Laplace's equation at every time step. Numerically solving Laplace's 
equation at every time step is expensive, particularly in three 
dimensions.

A characterization of $G(\eta)$ that avoids the expensive 
numerical solution of Lap\-lace's equation is appealing. In the current 
work we discuss four such characterizations developed in the 
context of water waves. These four methods are the operator expansion 
method of Craig \& Sulem (CS) \cite{CraigSulem}, the transformed
field expansion method (TFE) of Bruno \& Reitich \cite{brunoReitich}
and Nicholls \& Reitich \cite{NichollsReitich,NichollsReitich06},
the nonlocal implicit 
formulation of the DNO given by Ablowitz, Fokas \& Musslimani (AFM) 
\cite{AFM}, and a dual version to the AFM method, derived by Ablowitz \& 
Haut \cite{AblowitzHaut}, which we denote by AFM$^{*}$. Each of these 
methods has had remarkable theoretical utility in deriving reduced 
models for water waves in various physical regimes \cite{Craigetal}, 
in deriving conserved quantities \cite{AFM}, and also in providing the 
theoretical framework to pose some inverse problems 
\cite{OVDH,vasanDeconinck}. Additionally, each of these methods readily 
generalizes to the case of both varying bottom boundaries and three 
dimensional fluids.

Of course, many traditional methods to numerically solve Laplace's
equation exist, including the boundary integral method, conformal
mapping techniques, and the finite element method. In the present work
we chose to compare the CS and TFE operator expansions and the
AFM/AFM$^{*}$ nonlocal formulations with the boundary integral
method.  We do not consider conformal mapping \cite{dyachenko1996}, as
it does not extend to three dimensional problems, nor traditional
finite elements \cite{rycroft:13}, as the trade-off between high-order
elements and sparsity of the stiffness matrix makes them expensive
when high accuracy is desired.  Though it is not usually described in
this way, we regard the TFE method as a spectrally accurate variant of
the finite element method.  The boundary integral method is
particularly efficient in two dimensions for periodic problems as the
lattice sums involved have an explicit analytical
representation. However, in three dimensions, boundary integral
methods are considerably more difficult to implement.

The overall goal of the present work is to understand the relative
effectiveness and accuracy of each method for two-dimensional fluids.
The CS, AFM and AFM$^*$ methods involve highly ill-conditioned
intermediate calculations, so much of the paper focuses on whether
accurate results can be obtained if multiple-precision arithmetic is
employed in these intermediate calculations. In particular, for
timestepping the water wave, it is important that the methods work if
the input data is only known with limited accuracy --- additional
precision in intermediate calculations is permissible as long as 
the output is roughly as accurate as the input.  Given the care
required to obtain high accuracy in the present work, one must be
cautious about using these methods in double-precision without
carefully monitoring condition numbers and cancellation of digits.

An outline of the paper is as follows. In \S~\ref{DNO_reps}, we
briefly introduce the CS operator expansion, the AFM and AFM$^*$
nonlocal formulations, and the boundary integral and TFE
representations of the DNO. In \S\ref{assessment} we present
comparisons of the first four of these methods for specific choices of
the free surface $\eta$. Here we seek to quantify how accurate the CS,
AFM and AFM$^{*}$ methods are. We assume the exact DNO is obtained
from the boundary integral method, and, where suitable, perform
computations using higher precision. Section \ref{converge} discusses
the subtle cancellation properties associated with the series
representation of Craig \& Sulem. In \S\ref{examples}, we present a
specific example where the AFM$^{*}$ method (in its usual
interpretation as a system of equations for the coefficients of a
certain series) is guaranteed to fail, and other instances where AFM
and AFM$^{*}$ successfully converge to the correct Neumann data. The
subtle cancellation of the CS expansion is mirrored by the rapidly
decaying singular values associated with the AFM/AFM$^{*}$
methods. This leads us to consider
regularized versions of AFM/AFM$^*$ that involve oversampling the AFM
basis functions to accurately approximate a Gram-Schmidt
orthogonalization procedure via $QR$ factorization or the singular
value decomposition.  In the SVD approach, we also investigate the use
of a pseudo-inverse cutoff threshold.  We find that the AFM basis
functions can be more efficient at representing solutions than a
Fourier basis, but with the drawback of poor conditioning.  Finally,
in \S\ref{waterwave}, we discuss the behavior of all the methods on
examples relevant to water waves.

\section{Representations of the
  Dirichlet--Neumann operator}\label{DNO_reps}

As mentioned in the introduction, the Dirichlet--Neumann operator
plays an important role in the mathematical formulation of the motion
of surface gravity waves. To efficiently compute the time-dependent
motion, we require a fast and efficient means to solve Laplace's
equation, or alternatively, a direct method to compute the
Dirichlet--Neumann operator. In this section, we outline five commonly
used approaches. The first, due to Craig \& Sulem \cite{CraigSulem},
involves expanding the DNO in a Taylor series.  The second, due to
Ablowitz, Fokas and Musslimani \cite{AFM}, involves deriving a global
relation between the Dirichlet and Neumann data that can be used as an
integral equation to solve for the Neumann data.  The third
\cite{AblowitzHaut, vasanDeconinck, OliverasDeconinck, OVDH} is a dual
variant of the second, formulated more directly. The fourth is a
boundary integral collocation method \cite{lh76, baker:82, krasny:86,
  mercer:92, baker10, water1, water2}. And the fifth is the
transformed field expansion method of Bruno \& Reitich
\cite{brunoReitich} and Nicholls \& Reitich \cite{NichollsReitich,
NichollsReitich06}.

\subsection{Power Series Expansion of the DNO.}
\label{sec:cs:intro}

Consider Laplace's equation
\[\phi_{xx} + \phi_{yy} = 0,\]
posed on the domain $\Omega=\lbrace (x,y)\in\mathbb{R}^2:
0<x<L,-h<y<\eta(x)\rbrace$, where $\eta$ is a smooth periodic function
with period $L$. Further assume that $\phi$ is periodic in the
horizontal variable $x$ with period $L$, and
\[\phi_y(x,-h)=0.\]
Thus we restrict ourselves to a flat bottom boundary at $y=-h$.
%
%
Let $\mathscr D(x)$ and $\sN(x)$ be the Dirichlet and Neumann values of
the function $\phi$ at $y=\eta(x)$. If either $\sD(x)$ or $\sN(x)$ is given
(in appropriate function spaces), the problem of determining $\phi$ in
$\Omega$ is well-posed in the Hadamard sense.  As we require the map from
the Dirichlet to the Neumann data, assume we are given a Dirichlet
condition at the boundary $y=\eta(x)$. The associated Neumann
condition at the boundary $y=\eta(x)$ is given in terms of the
solution to the following boundary-value problem:
\begin{alignat*}{2}
\phi_{xx} + \phi_{yy} &= 0,& & -h<y<\eta(x), \\
\phi(x+L,y) &= \phi(x,y), \quad &&-h<y<\eta(x), \\
\eta(x+L)&=\eta(x),\\
\phi(x,\eta(x)) &= \sD(x), \\
\phi_y(x,-h) &= 0.
\end{alignat*}
In abstract terms, the Dirichlet--Neumann operator $G$ is given
by
\[  G(\eta)\mathscr{D} = \phi_{y} - \eta_{x}\phi_{x}, \]
where $\phi$ satisfies the above boundary-value problem. Note that a
function of the form
\[\varphi=\exp{(ikx)}\cosh(k(y+h)),\]
satisfies Laplace's equation, periodicity and the boundary condition
at $y=-h$ for $k=2\pi n/L,$ $n\in\mathbb{Z}$.  Hence
\begin{equation}\label{eq:G:hyp}
  G(\eta)\varphi(x,\eta) = ke^{ikx}\sinh(k(\eta+h)) -
  ik\eta_{x}e^{ik}\cosh(k(\eta+h)).
\end{equation}
It is well-known (see \cite{CraigSulem,NichollsReitich} and references
therein) that for a Lipschitz domain, the DNO is an analytic 
function of the domain shape. Thus, $G$ has a power series expansion 
in $\eta$.  Writing
\[G(\eta) = \sum_{j=0}^{\oo}G_j(\eta),\]
where $G_j(\lambda \eta) = \lambda^{j}G_j(\eta)$ for
$\lambda\in\mathbb{R}$, we obtain an explicit representation for
$G_j(\eta)$ from (\ref{eq:G:hyp}) by expanding the hyperbolic
terms in their respective Taylor series and identifying terms of the
same degree in $\eta$. To lowest order we obtain
\[G_0e^{ikx} = k \tanh(kh) e^{ikx}.\]
By decomposing the given Dirichlet condition in a Fourier series, we
obtain the following representation of the lowest order term of the
DNO
\[\mathcal{F}\big[ G_0\mathscr{D} \big] = k \tanh(kh)
\mathcal{F}\left[\mathscr{D}\right].\]
Similarly, proceeding to higher order, we obtain further terms in the
expansion of the DNO. For instance,
\begin{align*}
  G_1 &= D\eta D - G_0\eta G_0, \\
  G_2 &= -\frac{1}{2}\left(
    G_0\eta^2D^2 - 2G_0\eta G_0\eta G_0 +
    D^2\eta^2G_0\right),
\end{align*}
where $D=-i\partial_x$ and $G_0 = D\tanh(hD)$.  Terms such as
$\tanh(hD)$ are understood as pseudo differential operators,
\emph{i.e.}~they are defined through associated Fourier
multipliers. We remark that higher order terms of the Taylor expansion
of $G$ involve increasingly higher order derivatives. Although the
Taylor series of $G$ exists for $\eta$ with a limited (even finite)
degree of smoothness, the formulas for $G_j(\eta)$ are not valid in
such cases, or must be interpreted very carefully, as a whole, rather
than as a sum of individual operators. Even when $\eta$ is real
analytic, there is a delicate balance existing among the terms that
leads to a high degree of cancellation \cite{NichollsReitich}. This
causes numerical difficulties in finite precision arithmetic.  We
explore the extent of these cancellations in arbitrary precision
arithmetic in \S\ref{assessment}.  Of course, many of these problems
can be avoided by flattening the domain through a change of variables
\cite{brunoReitich,NichollsReitich,Lannes}, or using boundary integral
methods to compute the DNO \cite{water1,water2}. However, the original
Craig--Sulem expansion remains of theoretical interest, and is the one
most closely related to the AFM approach, discussed next.


\subsection{AFM Implicit Representation.}

Following \cite{AFM,AblowitzHaut}, we now derive the global relation
of the Ablowitz--Fokas--Musslimani reformulation of the water-wave
problem.  The distinguishing feature of this reformulation is the
implicit nonlocal characterization of the Dirichlet--Neumann operator.
Solving the resulting integral equation gives the full DNO,
effectively summing all the terms in the expansion of Craig \& Sulem
without having to compute them order by order.


As above, the functions
\begin{align}\label{psi}
    \psi=\exp{(ikx)}\cosh(k(y+h))
\end{align}
play a role, but now as dual functions rather than basis functions.
From Green's second identity, we have
\begin{align}\nonumber
0&=\int_D\left( \psi(\phi_{xx} + \phi_{yy}) - \phi(\psi_{xx} + \psi_{yy}) \right)dx\,dy,\\\nonumber
&= \int_{\partial D} \left(\psi\pp{\phi}{n} - \phi\pp{\psi}{n}\right) dS,\\\nonumber
&= \int_0^L \psi(x,\eta)\left[\phi_y(x,\eta)-\eta_x\phi_x(x,\eta)\right]dx -
\int_0^L\phi(x,\eta)\left[\psi_y(x,\eta)-\eta_x\psi_x(x,\eta)\right]dx,\\
&= \int_0^L\psi(x,\eta)\sN(x) dx - \int_0^L\sD(x)
\left[\psi_y(x,\eta)-\eta_x\psi_x(x,\eta)\right]dx,\label{Green}
\end{align}
where $\partial/\partial n$ is the normal derivative to the
surface. Using the definition of $\psi$ and noting that
\[e^{ikx}\left(k\sinh(k(\eta+h))-ik\eta_x\cosh(k(\eta+h))\right)
= -i\partial_x\left(e^{ikx}\sinh(k(\eta+h))\right),\]
we obtain
\begin{align}\label{AFM}
\int_0^Le^{ikx}\cosh(k(\eta+h)) \sN(x)\,dx =
\int_0^L ie^{ikx}\sinh(k(\eta+h))\partial_{x}\mathscr{D}\:(x)\,dx,
\end{align}
which is the Ablowitz-Fokas-Musslimani (AFM) global relation
\cite{AFM} for Laplace's equation. Note that the global relation thus
obtained is but a rephrasing of Green's identity. In \cite{AFM}, the
authors obtain this expression in a different, but equivalent, manner.
By a spectral collocation technique, this yields an algorithm in which
approximate values $\sN(x_j)$ are obtained on a grid
$\{x_j\}_{j=0}^{M-1}$ by solving a linear system. In that case,
(\ref{AFM}) is enforced for wave numbers $|k|\le (2\pi/L) M/2$.

A third approach, due to Ablowitz \& Haut \cite{AblowitzHaut}, can be
derived from (\ref{AFM}) through a type of inverse Fourier transform.
The resulting DNO algorithm boils down to solving
\begin{align}\label{Dirichlet}
  \sum_{k}\hat{\Psi}_{k}e^{ikx}\cosh(k(\eta+h)) =
  \sD(x),
\end{align}
for the expansion coefficients $\hat\Psi_k$, and then computing
\begin{align}\label{Normal}
  \sN(x) =
  -i\partial_x\left(\sum_{k}\hat{\Psi}_{k}e^{ikx}
    \sinh(k(\eta+h))\right).
\end{align}
%
%
For those functions $\sD(x)$ that permit an expansion of the form
(\ref{Dirichlet}), the normal derivative is expected (from term by
term differentiation) to be of the form (\ref{Normal}).

The operation of solving the integral equation (\ref{AFM}) is the
formal adjoint of solving the system (\ref{Dirichlet}),
(\ref{Normal}).  Indeed, if the (conjugate of the) former is written
$\sN=(A^*)^{-1}B^*\partial_x\sD$,
then the latter becomes
$\sN=-\partial_xBA^{-1}\sD$,
which are consistent since
$G(\eta)$ is self-adjoint.  Here
\begin{align}\notag
  Ac &= \sum_k c_k e^{ikx}\frac{\cosh(k(\eta(x)+h))}{w_k},
  \quad
  (A^*f)_k = \int_0^L e^{-ikx}\frac{\cosh(k(\eta(x)+h))}{w_k}f(x)\,dx,
  \\
  \label{wts}
  Bc &= \sum_k c_k ie^{ikx}\frac{\sinh(k(\eta(x)+h))}{w_k},
  \quad
  (B^*f)_k = \int_0^L e^{-ikx}\frac{\sinh(k(\eta(x)+h))}{iw_k}f(x)\,dx,
\end{align}
and the weights $w_k$ are chosen to make $A$ and $B$ bounded from
$\ell^2(\mathbb{Z})$ to $L^2(0,L)$.  Formally, these weights cancel
internally in the products $(A^*)^{-1}B^*$ and $BA^{-1}$.  Indeed, the
results are formally unchanged if both sides of (\ref{AFM}) are
multiplied by $w_k^{-1}$, or if the weights are absorbed into
$\hat\Psi_k$ in (\ref{Dirichlet}), (\ref{Normal}). These statements
are only formal since $A$ and $A^*$ are not invertible.

\subsection{Boundary Integral Method.}

Whereas the AFM$^*$ method represents $\phi$ in the fluid as a
superposition of basic solutions of the Laplace equation of the form
$e^{ky}e^{ikx}$, the boundary integral method represents $\phi$ as a
superposition of dipoles distributed along the surface,
\begin{equation*}
  \phi(z) = \int_{-\infty}^\infty -\der{N}{n_\zeta}(z,\zeta(\alpha))
  \mu(\zeta(\alpha))|\zeta'(\alpha)|\,d\alpha,
  \qquad N(z,\zeta) = \frac{1}{2\pi}\log|z-\zeta|.
\end{equation*}
Here $\zeta(\alpha) = \alpha + i\eta(\alpha)$ is a parametrization of
the free surface and $z=x+iy$ is a field point in the fluid. We then
use
\begin{equation*}
  -\der{N}{n_\zeta}\,ds = \opn{Im}\left\{\frac{\zeta'(\alpha)}{z - \zeta(\alpha)}
  \right\}d\alpha, \qquad
  \frac{1}{2}\cot\frac{z}{2} = PV\sum_k \frac{1}{z+2\pi k}
\end{equation*}
to reduce the integral to a period cell, and obtain
\begin{equation*}
  \phi(z) = \frac{1}{2\pi}\int_0^{2\pi}
  \widetilde{A}(z,\alpha)\mu(\alpha)\,d\alpha, \qquad
  \widetilde{A}(z,\alpha) = \opn{Im}\left\{
  \frac{\zeta'(\alpha)}{2}\cot\left(\frac{z-\zeta(\alpha)}{2}\right)\right\}.
\end{equation*}
Using the Plemelj formula \cite{muskhelishvili}, we take the limit
as the field point approaches the boundary from below to obtain a
second kind Fredholm integral equation for $\mu$:
\begin{align}
  \label{BIM1}
  &\frac{\mu(\alpha)}{2} + \frac{1}{2\pi}
  \int_0^{2\pi}A(\alpha,\beta)\mu(\beta)\,d\beta = \sD(\alpha), \\
  \notag
  &A(\alpha,\beta) = \opn{Im}\left\{
    \frac{\zeta'(\beta)}{2}
    \cot\left(\frac{\zeta(\alpha)-\zeta(\beta)}{2}\right) -
    \frac{1}{2}\cot\left(\frac{\alpha - \beta}{2}\right)
  \right\}.
\end{align}
Including $\frac{1}{2}\cot\left(\frac{\alpha - \beta}{2}\right)$ in
the formula has no effect on $A(\alpha,\beta)$, but shows that $A$ is
in fact a smooth function when $\eta$ is smooth.  Indeed, as
$\beta\rightarrow\alpha$, the $(\alpha-\beta)^{-1}$ singularities of
the terms in braces cancel, yielding $A(\alpha,\alpha)=\opn{Im}\big\{
-\zeta''(\alpha)/[2\zeta'(\alpha)]\big\}$.  Once $\mu$ is known,
the Neumann data is readily shown to satisfy
\begin{align}
  \label{BIM2}
  &\sN(\alpha) = \frac{1}{2}H[\mu'](\alpha)
  +\frac{1}{2\pi}\int_0^{2\pi}B(\alpha,\beta)\mu'(\beta)\,d\beta, \\
  \notag
  &B(\alpha,\beta) = \opn{Re}\left\{
    \frac{\zeta'(\alpha)}{2}
    \cot\left(\frac{\zeta(\alpha)-\zeta(\beta)}{2}\right) -
    \frac{1}{2}\cot\left(\frac{\alpha - \beta}{2}\right)
  \right\},
\end{align}
where $H$ is the Hilbert transform, with symbol $\hat H_k =
-i\opn{sgn}(k)$.  To carry this out numerically, $M$ collocation
points are used to turn the integral equation (\ref{BIM1}) into an
$M\times M$ matrix equation, where integrals are approximated by the
trapezoidal rule.  The derivative and Hilbert transform in
(\ref{BIM2}) are easily computed using the FFT.  The work involved in
setting up and solving these integral equations is very similar to
that of the AFM and AFM$^*$ methods.  However, the condition number is
much better in the BIM approach since the underlying infinite
dimensional system is a second-kind Fredholm integral equation.  This
makes a big difference in practice since intermediate calculations
need only be done in double-precision to achieve double-precision
results, and iterative methods such as GMRES can be employed to
reduce the work of solving the equations from $O(M^3)$ to $O(M^2)$.
See \cite{lh76, baker:82, krasny:86, mercer:92, mercer:94,
  baker:nachbin:98, smith:roberts:99, baker10, water2} for
similar boundary integral methods, including formulations that
incorporate a bottom boundary and allow the interface to overturn.

\subsection{Transformed Field Expansion method.}
\label{sec:tfe}

The aim of this approach
\cite{brunoReitich,NichollsReitich} is to compute successive terms
in the Craig--Sulem expansion via formulas that do not suffer from
catastrophic cancellation of digits in floating point arithmetic.  The
price we pay for this improvement is that the bulk fluid must be
discretized.  For simplicity, we consider only the finite depth case
in two dimensions. The three-dimensional case is considered in
\cite{NichollsReitich,NichollsReitich06}, while infinite depth is
treated in \cite{NichollsReitich06} by introducing a fictitious
interface coupling the unbounded problem on a half-space to the
finite-depth problem with a curved upper boundary and a flat lower
boundary.

Instead of deriving the perturbation expansion for $G(\eta)$ using
the ill-conditioned basis functions $e^{ikx}e^{\pm ky}$, as was done
in (\ref{eq:G:hyp}) above, we perform a boundary-flattening
change of variables:
\begin{equation*}
  u(x,y) = \phi\big(x,(1+h^{-1}\eta)y+\eta\big), \qquad
  0\le x<L, \qquad -h<y<0.
\end{equation*}
A straightforward calculation reveals that
\begin{align*}
  u_{xx} &= \phi_{xx} + \big(1+h^{-1}y\big)\big(1+h^{-1}\eta\big)^{-1}\eta_x \partial_y
  [ u_x - \circledast ] + \partial_x[ \circledast ], \\
  u_{yy} &= \phi_{yy} + \Big[1 - \big(1 + h^{-1}\eta \big)^{-2}\Big] u_{yy},
\end{align*}
where 
$\circledast = (1+h^{-1}y)(1+h^{-1}\eta)^{-1}\eta_x u_y$.
Since $\Delta\phi=0$, we find that
\begin{equation*}
  \Delta u = \partial_x\big( F_1 \big) + \partial_y\big( F_2 \big) + F_3,
\end{equation*}
where
\begin{align*}
  F_1 &= \big(1 + h^{-1}y\big)E\eta_xu_y, \qquad   E = \big(1+h^{-1}\eta\big)^{-1}, \\
  F_2 &= \big(1 + h^{-1}y\big)E\eta_xu_x -
  \big(1 + h^{-1}y\big)^2E^2\eta_x^2u_y + \big(1-E^2\big)u_y, \\
  F_3 &= -h^{-1}E\eta_xu_x + h^{-1}(1+h^{-1}y)E^2\eta_x^2u_y.
\end{align*}
Next, we write $\eta(x)=\veps f(x)$ and expand $E$ and $E^2$ in powers of
$\veps$ to conclude that the terms of the series $u(x,y) =
\sum_{n=0}^\infty \veps^n u_n(x,y)$ satisfy
\begin{alignat}{2} \label{lap:u0}
  \Delta u_0 &= 0, &
  u_0(x,0) &= \sD(x), \\ \label{lap:un}
  \Delta u_n &= \partial_x \big(F_1^n\big) + \partial_y\big(F_2^n\big) + F_3^n, 
  \qquad &
  u_n(x,0) &= 0,
\end{alignat}
as well as $u_{n,y}(x,-h)=0$ and $u_n(x+L,y) = u_n(x,y)$.  Here
\begin{align*}
  F_1^n &= (1 + h^{-1}y)f_x \sum_{m=0}^{n-1}(-h^{-1}f)^m
  u_{n-1-m,y}, \\
  F_2^n &= (1+h^{-1}y)(F_4 - F_5) + \frac{f}{h}
  \sum_{m=0}^{n-1}(m+2)(-h^{-1}f)^mu_{n-1-m,y}, \\
  F_3^n &= h^{-1}(F_5 - F_4), \qquad
  F_4^n = f_x\sum_{m=0}^{n-1}(-h^{-1}f)^mu_{n-1-m,x}, \\
  F_5^n &= (1+h^{-1}y)f_x^2\sum_{m=0}^{n-2} (l+1)(-h^{-1}f)^m u_{n-2-m,y}.
\end{align*}
Finally, we use
%
\begin{align*}
  &G(\eta)\sD = \mathbf{n}\cdot\big[\nabla\phi\big]_{y=0} =
  \Big[-\eta_xu_x + \big(1 + \eta_x^2)(1 + h^{-1}\eta)^{-1}u_y\Big]_{y=0}, \\
  &\mathbf{n} = (-\eta_x,1), \qquad
  \nabla\phi = \Big( u_x - (1+h^{-1}y)(1+h^{-1}\eta)^{-1}\eta_x u_y\,,\,
  (1+h^{-1}\eta)^{-1}u_y \Big)
\end{align*}
to conclude
that $G(\veps f) = \sum_{n=0}^\infty \veps^n G_n(f)$ with
\begin{equation*}
  G_n(f)\sD = -f_x u_{n-1,x} +
  \sum_{m=0}^n(-h^{-1}f)^m u_{n-m,y} + f_x^2\sum_{m=0}^{n-2}
  (-h^{-1}f)^m u_{n-2-m,y},
\end{equation*}
where empty sums (with upper index smaller than lower index) are zero.

In our code, $u_{n,x}(x,y)$ and $u_{n,y}(x,y)$ are stored on a
rectilinear grid with $M$ uniformly spaced mesh points in the
$x$-direction and $N+1$ Chebyshev-Lobatto nodes in the $y$-direction,
mapped by an affine transformation to obtain $y_0=-h$ and $y_N=0$.
Functions on the grid are stored as matrices with entries in a column
indexed by $x$, holding $y$ fixed.  The formulas for $F^n_j(x,y)$ are
evaluated pointwise on the grid from the known values of $f(x)$,
$f_x(x)$, $u_{m,x}(x,y)$ and $u_{m,y}(x,y)$ for $m=0,\dots,n-1$.  The
FFT of the zeroth order term $u_0(x,y)$ is computed from (\ref{lap:u0}) by
expanding $\sD(x)$ in a Fourier series and evaluating $\hat u_0(k,y) =
\hat\sD_k \cosh(k(y+h))\opn{sech}(kh)$ at the grid points $y_j$.  To
obtain $u_{0,x}(x,y)$, we multiply $\hat u_0(k,y)$ by $ik$ and take the
inverse FFT.  Similarly, $u_{0,y}(x,y)$ is obtained by taking the iFFT
of $\partial_y\hat u_0(k,y)$. This latter function is computed by
transforming the rows of $\hat u_0$ to their Chebyshev coefficients
(also using an FFT), differentiating the Chebyshev polynomials, and
evaluating the result on the grid using the Clenshaw recurrence
formula. The differentiation procedure amounts to determining the
coefficients of $\partial_y\hat u_{0}$ from those of $\hat u_0$:
%
\begin{align*}
  &\hat u_0(k,y)=\jt\sum_j\alpha_j(k)T_j(1+2h^{-1}y), \quad
  \partial_y \hat u_0(k,y) = \sum_j\beta_j(k)T_j(1+2h^{-1}y)2h^{-1}, \\
  &\beta_N=0, \quad \beta_{N-1}=N\alpha_N, \quad
  \beta_j = (j+1)\alpha_{j+1} + \beta_{j+2}, \quad (j=N-2:-1:0),
\end{align*}
where $T_j(y)$ is the $j$th Chebyshev polynomial.
Finally, (\ref{lap:un}) is solved using an FFT in
the $x$-direction to convert the PDE into a collection of uncoupled
boundary value problems in $y$. The $k$th BVP is
\begin{equation*}
  \partial_y^2 \hat u_n(k,y) - k^2 \hat u_n(k,y) =
  ik\hat F_1^n(k,y) + \partial_y \hat F_2^n(k,y) + \hat F_3^n(k,y),
\end{equation*}
subject to $\partial_y\hat u_n(k,-h)=0$, $\hat u_n(k,0)=0$.  Instead
of using the Chebyshev tau method \cite{NichollsReitich06, canuto:88}
to solve this BVP, we multiply by a test function, integrate the first
term by parts, and proceed as if implementing a finite element method
using Chebyshev polynomials as the basis functions. Our implementation
is similar to what was done in \cite{plasma2} to study the projected
dynamics of kinetic diffusion equations in spaces of orthogonal
polynomials.  Once $\hat u_n(k,y)$ is known for each $k$, we compute
$u_{n,x}$ and $u_{n,y}$ as described above for $u_{0,x}$ and
$u_{0,y}$.

The key observation is that the derivatives on $F_1^n$ and $F_2^n$ in
(\ref{lap:un}), and on $u_n$ at the end of the procedure (to obtain
$u_{n,x}$ and $u_{n,y}$), are balanced by the inverse Laplacian in
(\ref{lap:un}). In other words, at the point in the algorithm where
large numbers enter the computation due to applying derivatives in
Fourier space, we divide by even larger numbers by applying the
inverse Laplacian.  By contrast, in the CS expansion, $G_n(f)$ is
expressed as a sum of several terms, each involving $n$ derivatives of
products of $\sD$ with powers of $\eta$. Each derivative amplifies
roundoff error, making it difficult to extract the desired sum, which
is often many orders of magnitude smaller than the individual terms,
as shown below.

\section{Comparison of the methods.}
\label{assessment}

In this section we assess the merits and shortcomings of the CS method
and the AFM/AFM$^*$ methods.  In \S\ref{converge}, we explore the
convergence of the CS expansion and quantify the delicate cancellation
of terms mentioned in \cite{NichollsReitich}.  In \S\ref{examples}, we
present an example illustrating a circumstance in which the AFM$^{*}$
method (in its usual interpretation) is guaranteed to fail. By
oversampling the columns of the linear operators involved in the
AFM/AFM$^{*}$ methods, we are once again able to compute the true
normal derivative with spectral accuracy.  We also report on the
performance of the CS expansion method for the examples considered for
the AFM method. Finally in \S\ref{waterwave} we consider the
performance of these three methods, as well as the TFE and BIM
methods, on two examples in which the Dirichlet data and free surface
come from solutions of the water-wave equations. The first is a
standing water wave and the second involves two interacting traveling
water waves.

\subsection{Cancellation properties of the DNO expansion.}
\label{converge}

For simplicity, until \S\ref{waterwave}, we restrict to the case of
two-dimensional fluids of infinite depth.  In that case, the DNO
expansion takes the form \cite{NichollsReitich}
\begin{align}
  G_0(f) &= |D|, \\
  G_n(f) &= |D|^{n-1}D\frac{f^n}{n!}D - \sum_{s=0}^{n-1} |D|^{n-s}\frac{f^{n-s}}{(n-s)!}G_s(f), \quad
  n=1,2,3,\dots,
\end{align}
where the symbols of $D$ and $|D|$ are $k$ and $|k|$, respectively.  
Note that $G_n(\eta) = \veps^n G_n(f)$ when $\eta = \veps f$.
Computationally, it is convenient to absorb the $s=0$ term into
the first term:
\begin{equation}\label{recur}
  G_n(f) = A_n(f) - \sum_{s=1}^{n-1} \frac{1}{(n-s)!}|D|^{n-s} f^{n-s} G_s(f),
  \quad n=1,2,3,\dots,
\end{equation}
where $A_n(f) =\frac{1}{n!} |D|^{n-1}\big( Df^n D - |D|f^n |D|\big)$ for $n\ge1$.
In Fourier space, $A_n(f)$ is an infinite matrix with the quadrants containing
the main diagonal zeroed out:
\begin{equation}\label{eq:An:formula}
  A_n(f)^\wedge_{kj} = \begin{cases} \jd
    -\frac{2|k|^n|j|}{n!}(f^n)^\wedge_{k-j}, & kj<0, \\
    0 & kj\ge0.
  \end{cases}
\end{equation}
This already accomplishes a fair amount of cancellation since
$|k-j|>|k|$ when $jk<0$, so the rapid growth of $|k|^n/n!$ is
balanced by decay of $(f^n)^\wedge_{k-j}$ when $f$ is smooth.
Indeed, if $f$ is real analytic, one may show that there exist
$C$ and $\rho$ such that $|(f^n)^\wedge_k|\le C^ne^{-\rho|k|}$,
which is enough to guarantee that $A_n(f)^\wedge$ maps $l^2$
sequences $\hat\varphi$ to exponentially decaying sequences:
\begin{equation*}
  \big|\big( A_n(f)\varphi \big)^\wedge_k\big| \le
  \frac{|Ck|^ne^{-\rho|k|}}{n!}
  \sum_{j\in J(k)} 2|j|e^{-\rho|j|}\,|\hat\varphi_j| \le
  \left(\rho^{-3/2}e^{2\rho}\|\hat\varphi\|_{l^2}\right)
  \frac{|Ck|^n}{n!}e^{-\rho|k|}.
\end{equation*}
Here $J(k)$ is the set of positive integers when $k$ is
negative and the set of negative integers when $k$ is positive.
Thus, if $f$ is real analytic, (\ref{recur}) implies that
$G_n(f)\varphi$ is real analytic for $n\ge1$ when $\varphi$
is merely $L^2$.  Also, aside from $n=0$, $G_n(f)$ is bounded
on $L^2$ when $f$ is real analytic.

\begin{figure}
    \includegraphics[width=\linewidth,trim=30 0 10 -10]{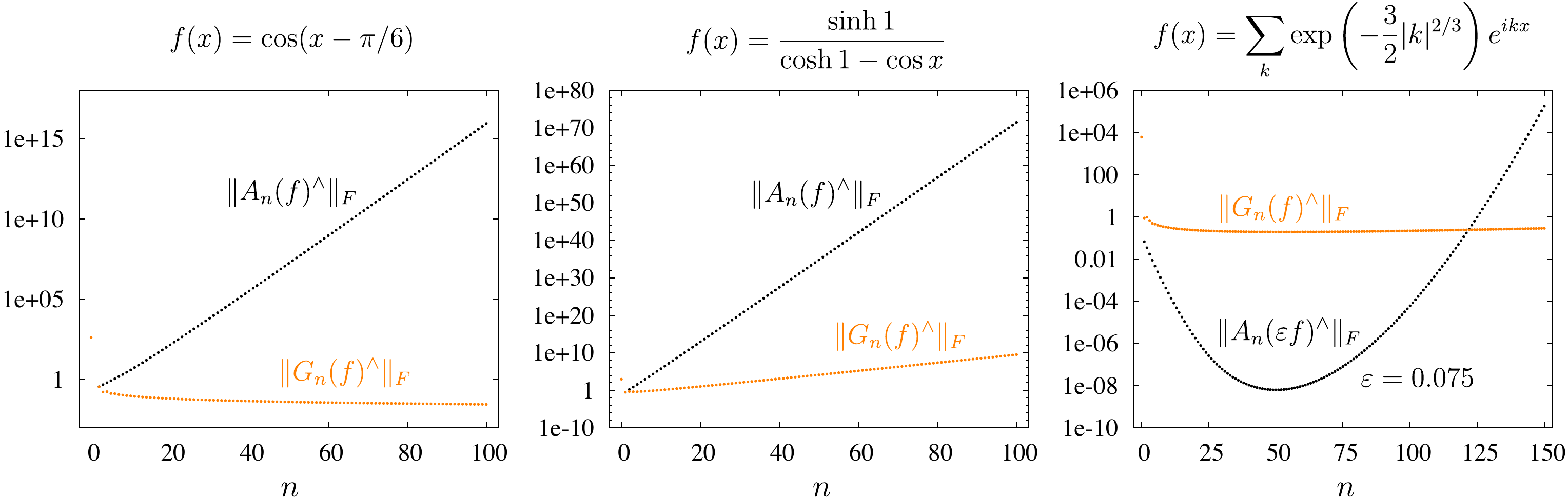}
    \caption{\label{fig:growth} Illustration of the cancellations that
      occur in the recursion (\ref{recur}) for three types of
      functions: band-limited, real-analytic, and $C^\infty$. In the
      third panel, the Frobenius norm of $A_n(f)^\wedge$ grows so much
      faster than that of $G_n(f)^\wedge$ that we had to
      re-scale it by $\veps^n$, $\veps=0.075$.  Note that
      $\veps^{150}\approx1.8\times10^{-169}$.}
\end{figure}

Numerical experiments reveal that significant additional cancellations
occur in (\ref{recur}), beyond combining $Df^nD$ with $|D|f^n|D|$.
Indeed, without these cancellations, even when $f$ is real analytic,
one would expect $(|D|f)^{n-1}A_1(f)$, which is one of the terms in
(\ref{recur}) when the recursion is unrolled, to grow
super-exponentially with $n$. In Figure~\ref{fig:growth}, we plot
the Frobenius norm of the operators $A_n(f)^\wedge$ and $G_n(f)^\wedge$
versus $n$ for the functions
\begin{align}
  \text{Example 1:}& \qquad f(x)=\cos(x-\pi/6), \\
  \label{ex2}
  \text{Example 2:}& \qquad f(x)=\frac{\sinh(1)}{\cosh(1)-\cos(x)}, \\
  \label{ex3}
  \text{Example 3:}& \qquad
  f(x)=\sum_k \exp\left(-\frac{3}{2}|k|^{2/3}\right)e^{ikx},
\end{align}
which are band-limited, real-analytic, and $C^\infty$, respectively.
The Frobenius norm of a matrix (in this case doubly-infinite) is the
root sum of squares of the matrix entries. We computed it in arbitrary
precision arithmetic using (\ref{eq:An:formula}) to evaluate
$A_n(f)^\wedge$ and (\ref{recur}) to evaluate $G_n(f)$. The computations
were done column by column, in parallel.  We worked in Fourier space
except when computing $f^{n-s}G_s(f)$, which was done by applying the
inverse FFT to a column of $G_s(f)^\wedge$, then multiplying by
$f^{n-s}$ in real space, and finally applying the FFT again. The
Frobenius norm of $A(f)^\wedge$ was computed from the indices in
the range
\begin{equation}\label{finite:lattice}
  k<0, \qquad j>0, \qquad |k-j|<M/2,
\end{equation}
where $M$ was chosen large enough that, for the range of $n$'s
considered, the terms $A_n(f)^\wedge_{kj}$ with $|k-j|\ge M/2$ are
small compared to the largest of those satisfying $|k-j|<M/2$, and may
be set to zero. As explained above, this is possible since the
exponential decay of $(f^n)^\wedge_{k-j}$ dominates the polynomial
growth of $|k|^n|j|$ in (\ref{eq:An:formula}).  This $M$ was also used
as the number of grid points in the FFT. We include a factor of 2 when
summing the squares of the matrix entries to account for the entries
$A_n(f)^\wedge_{-k,-j} = \overline{A_n(f)^\wedge_{kj}}$ in the
opposite quadrant, $k>0$, $j<0$.

The other consideration for choosing $M$ is that errors near the
boundary propagate inward when computing $G_n(f)^\wedge_{kj}$.  Thus,
we choose a smaller integer $K$ and compute the Frobenius norm of
$G_n(f)^\wedge$ from the entries with indices
\begin{equation}\label{eq:G:range}
  -K/2<k<K/2, \qquad 0<j<K/2.
\end{equation}
The remaining columns of $G_n(f)^\wedge$ (with $K/2\le j<M/2$) are not
computed, although the rows are computed out to $-M/2<k<M/2$.  We
always zero out the Nyquist frequency, $|k|=M/2$.  As before,
$G_n(f)^\wedge_{-k,-j} = \overline{G_n(f)^\wedge_{kj}}$ is accounted
for in Figure~\ref{fig:growth} with a factor of $\sqrt{2}$ in the root
sum of squares.

\begin{figure}
    \includegraphics[width=\linewidth,trim=30 0 10 -10]{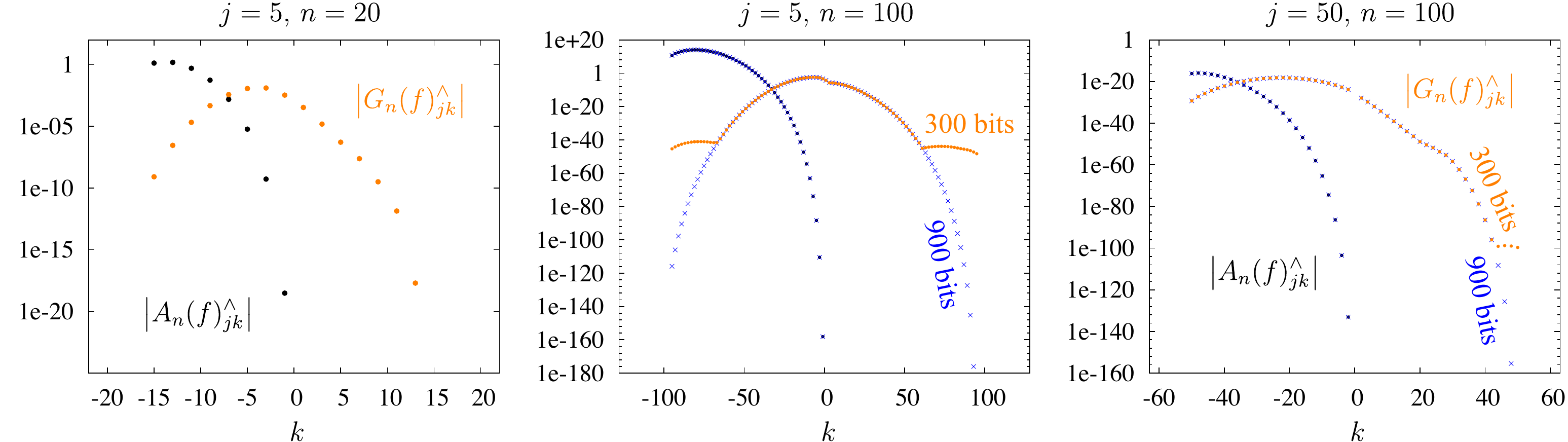}
    \caption{\label{fig:cosA} Plots of the magnitudes of the non-zero
      matrix entries in selected columns (indexed by $j$) of
      $A_n(f)^\wedge$ and $G_n(f)^\wedge$ for Example 1.  The orange
      and black markers were computed with 300 bits of precision while
      the blue markers were computed with 900 bits of precision.  The
      support of each column is finite since $(f^n)^\wedge_k=0$ for
      $|k|>n$.}
\end{figure}

The parameters used in these computations were
\begin{equation*}
  \begin{array}{c|c|c|c|c}
    \text{Example} & M & K & n_\text{max} & \text{bits} \\ \hline
    1 & 256 & 128 & 100 & 300 \\
    2 & 2048 & 330 & 100 & 500 \\
    3 & 24576 & 768 & 150 & 1500
  \end{array}
\end{equation*}
Here ``bits'' refers to the binary precision of the mantissa, where 53
would correspond to double-precision.  We used MPFR \cite{mpfr} for
the floating-point arithmetic, which provides IEEE-like arbitrary
precision rounding behavior.  We also used double-precision and
quadruple-precision arithmetic (using the qd package) in some cases.
The following table gives running times for multiplying two
$1000\times1000$ matrices on a 3.33 GHz Intel Xeon X5680 system with
12 cores:
\begin{equation*}
  \begin{array}{r|c|c|c|c|c}
    \text{precision} & \text{double} & \text{quad} & \text{300 bits} &
    \text{500 bits} & \text{1500 bits} \\ \hline
    \text{time (in seconds)} & 0.0157 & 1.66 & 12.5 & 17.1 & 63.5
  \end{array}
\end{equation*}
The double-precision calculation is particularly fast due to the use
of Intel's math kernel library. The higher-precision examples were
parallelized using openMP, but do not employ block-matrix algorithms
to re-use data that has been pulled from main memory to cache. Our
general experience (excluding level 3 BLAS routines such as
matrix-matrix multiplication) is that switching from double to
quadruple to arbitrary precision slows down the calculation by a
factor of 10 each.

\begin{figure}
  \includegraphics[width=\linewidth,trim=30 0 10 -10]{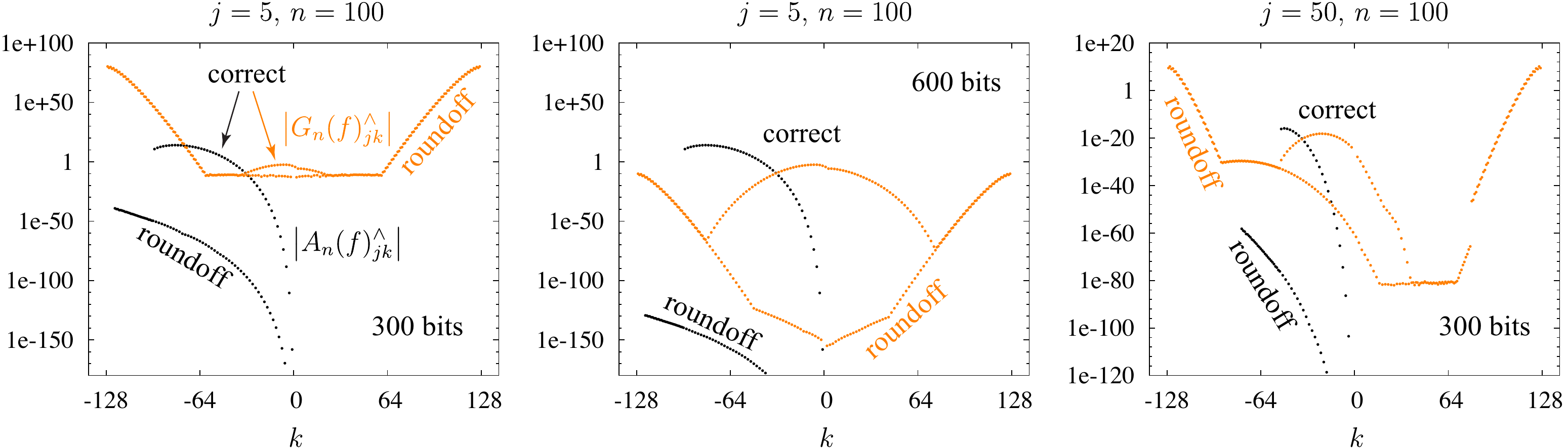}
  \caption{\label{fig:cosB} Repeat of the calculation of
    Figure~\ref{fig:cosA} without imposing a zero-pattern filter on
    the matrices as they are constructed. The $n=20$ solution has been
    replaced by a second instance of the $n=100$ solution.  (left and
    right) With 300 bits of precision, $G_n(f)^\wedge$ is almost
    entirely corrupted with roundoff errors. (center) With 600 bits of
    precision, roundoff error is suppressed enough to achieve an
    accurate result.  Errors are largest near $k=\pm M/2$, and
    propagates inward as $n$ increases.}
\end{figure}

In Figure~\ref{fig:cosA}, we plot the
non-zero matrix entries of $A_n(f)^\wedge$ and $G_n(f)^\wedge$ for
Example 1.  Because $f(x)=\cos(x-\pi/6)$, the Fourier modes
$(f^n)^\wedge_k$ are zero for $|k|>n$ or $k-n$ odd.  As a result,
$A_n(f)^\wedge$ and $G_n(f)^\wedge$ have only finitely many nonzero
terms in this example:
\begin{alignat*}{3}
  A_n(f)^\wedge_{kj}&=0 \quad & &\text{if} &\quad kj\ge0, \;\; &|k|>n-|j|, \;\; \text{or} \;\;
  k-j-n \;\; \text{is odd}, \\
  G_n(f)^\wedge_{kj}&=0 \quad & &\text{if} &\quad kj=0, \;\; &|k|>n-|j|, \;\; \text{or} \;\; k-j-n \;\;
  \text{is odd}.
\end{alignat*}
%
In Figure~\ref{fig:cosA}, we explicitly filtered the data to zero out
matrix entries of $A_n(f)^\wedge$ and $G_n(f)^\wedge$ with indices in
these ranges.  If this is not done, roundoff error from the FFT is
rapidly amplified by the recurrence (\ref{recur}), and requires
additional precision to maintain accuracy.  This is demonstrated in
Figure~\ref{fig:cosB}, where we did not filter the data. Increasing
the precision from 300 bits to 600 causes the correct values of
$G_n(f)^\wedge$ to emerge from the roundoff noise.  This is less
of an issue for $A_n(f)^\wedge$, which involves errors from taking
the FFT of $f^n$, amplified by $|k|^n/n!$, but no recurrence.

\begin{figure}
  \includegraphics[width=\linewidth,trim=30 0 10 -10]{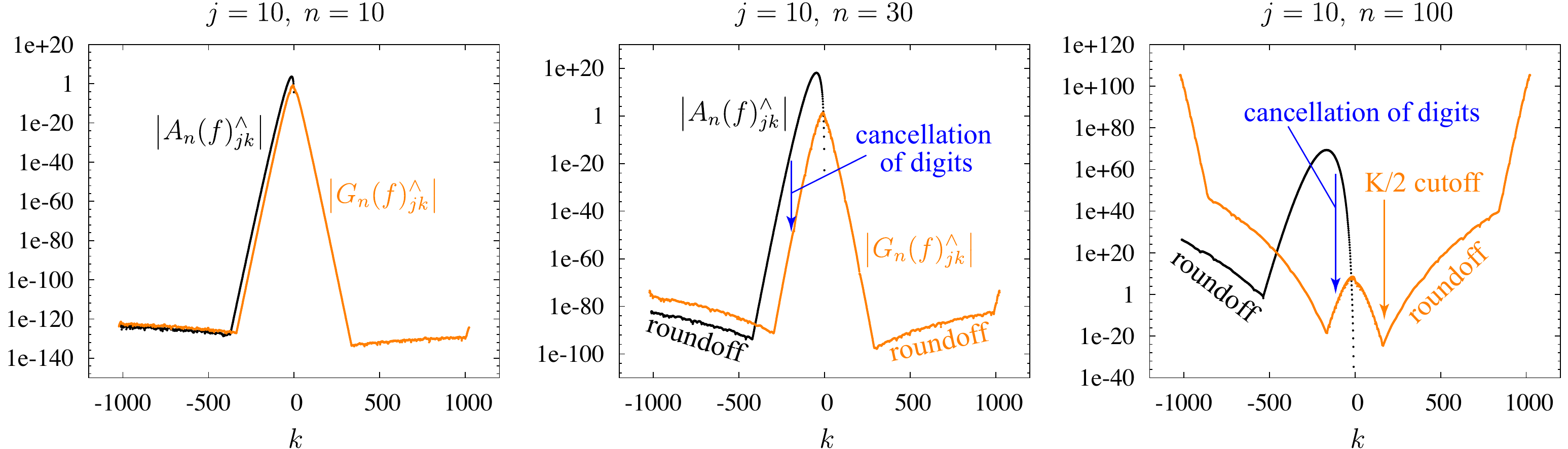}
  \caption{\label{fig:alpha} Plots of the magnitudes of the matrix
    entries in the $j=10$th column of $A_n(f)^\wedge$ and
    $G_n(f)^\wedge$, with $n=10,30,100$, for Example 2. For small $n$,
    the amplitude of $G_n(f)^\wedge_{kj}$ is similar to that of
    $A_n(f)^\wedge_{kj}$, indicating that little cancellation has
    occurred.  By the time $n$ reaches 100, the leading 90 digits of
    $A_n(f)^\wedge_{kj}$ have been eliminated in the recurrence
    (\ref{recur}) to obtain $G_n(f)^\wedge_{kj}$ for typical values
    of $j$, $k$. }
\end{figure}

\begin{figure}
  \includegraphics[width=\linewidth,trim=30 0 10 -10]{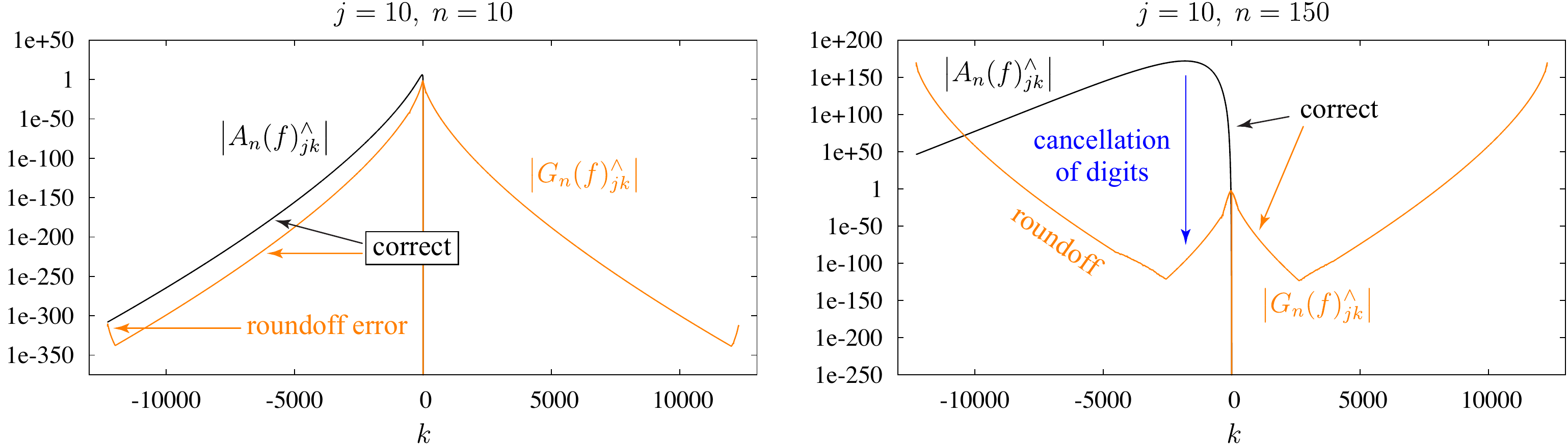}
  \caption{\label{fig:two3rds} Plots of the magnitudes of the matrix
    entries in the $10$th column of $A_n(f)^\wedge$ and
    $G_n(f)^\wedge$, with $n=10$ and $150$, for Example 3. $M$ had to
    be chosen quite large ($M=24576$) to prevent errors from
    propagating inward to $k=0$ before $n=150$. Note that more than
    250 leading digits of $A_n(f)^\wedge$ cancel to form
    $G_n(f)^\wedge$ for this column when $n=150$.  }
\end{figure}

In Figures~\ref{fig:alpha} and~\ref{fig:two3rds}, we plot selected
columns of $A_n(f)^\wedge$ and $G_n(f)^\wedge$ for Examples~2 and~3.
The main change from the band-limited case of Example~1 is that the
matrices are no longer of finite rank, and do not have compactly
supported columns. For small $n$, the decay rate of each column (with
respect to row index $k$) is still very fast, and the cancellations in
obtaining $G_n(f)^\wedge$ from $A_n(f)^\wedge$ is fairly mild.
However, as $n$ increases, the cancellations become quite severe.  The
cancellations can be seen in the figures as the vertical difference
from one curve to the other.  Recall that $A_n(f)^\wedge$ is the first
term in the formula (\ref{recur}) for $G_n(f)^\wedge$.  The functions
(\ref{ex2}) and (\ref{ex3}) have Fourier modes of the form
\begin{equation}\label{eq:fab}
  \text{Examples 2 and 3:} \qquad \hat f_k = e^{-\alpha |k|^\beta},
\end{equation}
where $\alpha=\beta=1$ in Example 2 and $\alpha=\beta^{-1}=3/2$ in
Example 3.  Example 2 is intended to represent a typical real-analytic
function, while Example 3 was designed to check if super-exponential
growth in the norms of the operators $A_n(f)$ might cause $G_n(f)$ to
also grow super-exponentially with $n$.  To see why $A_n(f)$ grows
super-exponentially, note that iterated convolution of $\hat f$ with
itself yields $(f^n)^\wedge_k\ge\hat f_k = e^{-\alpha |k|^\beta}$ for
functions of the form (\ref{eq:fab}).  We then consider the $j=1$
column of (\ref{eq:An:formula}) and maximize
\begin{equation*}
  \max_k A_n(f)^\wedge_{k,1} \ge \max_k |k|^n e^{-\alpha |k|^\beta}/n!.
\end{equation*}
The maximum on the right occurs near $k^*=-(n/\alpha\beta)^{1/\beta}$,
so we set $\alpha=\beta^{-1}$ and obtain, via Sterling's formula,
\begin{equation*}
  A_n(f)^\wedge_{k^*,1} \ge (\sqrt{2\pi n})^{\beta^{-1}}(n!)^{\beta^{-1}-1}.
\end{equation*}
We tried $\beta=1/2$ and $\beta=2/3$.  The decay in the former case
was too slow for the problem to be computationally tractable beyond
$n=50$.  So we present the results in Figure~\ref{fig:two3rds} with
$\beta=2/3$, where we were able to compute terms out to $n=150$.
Returning to Figure~\ref{fig:growth}, the right panel shows that the
Frobenius norm of $A_n(f)^\wedge$ does indeed grow
super-exponentially, but the cancellations are strong enough that
$G_n(f)^\wedge$ remains fairly flat.  As a result, $\|G_n(\veps
f)^\wedge\|_F$ decays like $\veps^n$ for large $n$ even though
$\|A_n(\veps f)^\wedge\|_F$ eventually stops decaying for any positive
$\veps$.

\begin{figure}
  \includegraphics[width=\linewidth,trim=30 0 10 -10]{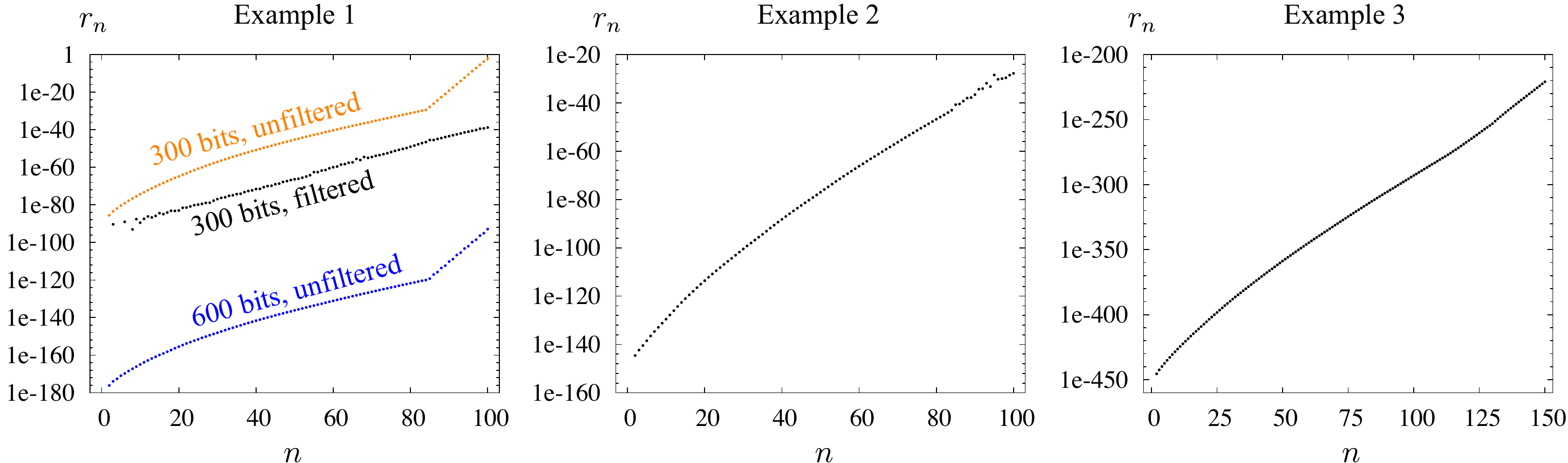}
  \caption{\label{fig:sym} Plot of symmetry errors in the matrix
    $G_n(f)^\wedge$ versus $n$. While filtering the data improves the
    band-limited case significantly, the more important factor is the
    precision of the underlying floating point arithmetic. }
\end{figure}

In addition to visual confirmation that roundoff error has not
corrupted $G_n(f)^\wedge_{jk}$ for $|k|$ small, as demonstrated in
Figures~\ref{fig:cosB}--\ref{fig:two3rds}, we also validate the
results by checking self-adjointness of each $G_n(f)^\wedge$.  This is
done by measuring
\begin{equation}
  r_n = \frac{\max_{j,k\in \mathbb{K}}
  \left|G_n(f)^\wedge_{kj} - \overline{G_n(f)^\wedge_{jk}}\right|}{\|G_n(f)^\wedge\|_F},
  \qquad  \mathbb{K} = \{k\;:\;|k|<K/2\}.
\end{equation}
In Figure~\ref{fig:sym}, we plot $r_n$ versus $n$ for each of the
three examples.  Note that 300 bits of precision was not sufficient
in the unfiltered case of Example 1 to avoid $O(1)$ errors in the
symmetry of $G_n(f)^\wedge$, whereas 600 bits gives at least 90
correct digits.  This is consistent with the results of Figure~\ref{fig:cosB},
where the signal is barely distinguishable from the noise in the left
panel, but is many orders of magnitude larger in the center panel.
In hindsight, 1500 bits was overkill for Example 3 since the symmetry
errors in $G_n(f)^\wedge$ are still below $10^{-200}$ when $n=150$.
Further validation of the correctness of the DNO expansion will be
given in the following section, where it will be used to compute
$G(\eta)\sD$ in a case where the solution is known.

\subsection{Examples illustrating success and failure of the AFM method.}
\label{examples}

In this section we present examples which show that even for
band-limited wave profiles $\eta$ of arbitrarily small amplitude, the
system (\ref{Dirichlet}-\ref{Normal}), with the sums interpreted as
infinite series, may fail to produce the Neumann data for Dirichlet
data of a certain class, which includes real analytic functions,
regardless of how many digits of arithmetic are used in the
computation.

Both AFM and AFM$^{*}$ methods involve linear operators acting on the
unknown Neumann data. In this section, we study the singular value
decomposition (SVD) of truncations of these linear operators. By
over-sampling the columns, we obtain better approximations of the SVD
of the underlying quasi-matrices $A$ and $B$ in (\ref{wts}), whose
columns are continuous functions.  This leads to regularized versions
of the AFM/AFM$^{*}$ methods, the simplest version (with no
pseudo-inverse cutoff) being equivalent to performing a Gram-Schmidt
orthogonalization of the basis functions $e^{ikx}\cosh(k(\eta(x)+h))$
before attempting to represent $\sD(x)$ or $\sN(x)$ using these
functions.  The regularized approach enables these methods to be used
even if no series of the form (\ref{Dirichlet}) exists for $\sD(x)$.

We now construct a function $\sD(x)$ for which no such series exists.
Consider the function 
\begin{align}
\label{wilken_ex}
\phi(x,y) &= \frac 1 2\mathrm{Im}\left\lbrace \cot\left(
    \frac{x+iy}{2}\right) -
  \cot\left(\frac{x+iy+2ih}{2}\right) \right\rbrace,
\end{align}
%
which is $2\pi$-periodic in $x$ and harmonic for all
$x,y$ outside of the set $(2\pi\mathbb{Z})\times\{0,-2h\}$.
Further,
$\phi(x,-h-y)=\phi(x,-h+y)$ and hence $\phi_y(x,-h)=0$ for all
$x$. Evaluating at $y=-h$, we obtain
\begin{align*}
  \phi(x,-h) &= \frac{\sinh(h)}{\cosh(h)-\cos(x)}
  = 1 + 2\sum_{k=1}^\infty e^{-kh}\cos kx =
  \sum_{k=-\infty}^\infty e^{-|k|h}e^{ikx}.
\end{align*}
For values of $y$ in the range $-2h<y<0$, we may also write
\begin{equation}\label{phi:ex1}
  \begin{aligned}
  \phi(x,y) &= 1 + 2\sum_{k=1}^\infty e^{-kh}\cosh(ky+kh)\cos kx \\
  &= 1 + \sum_{k=1}^\infty (e^{ky} + e^{-k(y+2h)})\cos kx, \qquad (-2h<y<0).
  \end{aligned}
\end{equation}
In the form (\ref{wilken_ex}), $\phi$ is well-behaved except at the
poles.  However, the AFM formulation is based on representing $\sD$
via the series (\ref{phi:ex1}), which is divergent for $y\ge0$.  Thus,
we expect trouble for wave profiles $\eta(x)$ that extend above
$y=0$.

\begin{figure}
\includegraphics[width=.9\linewidth]{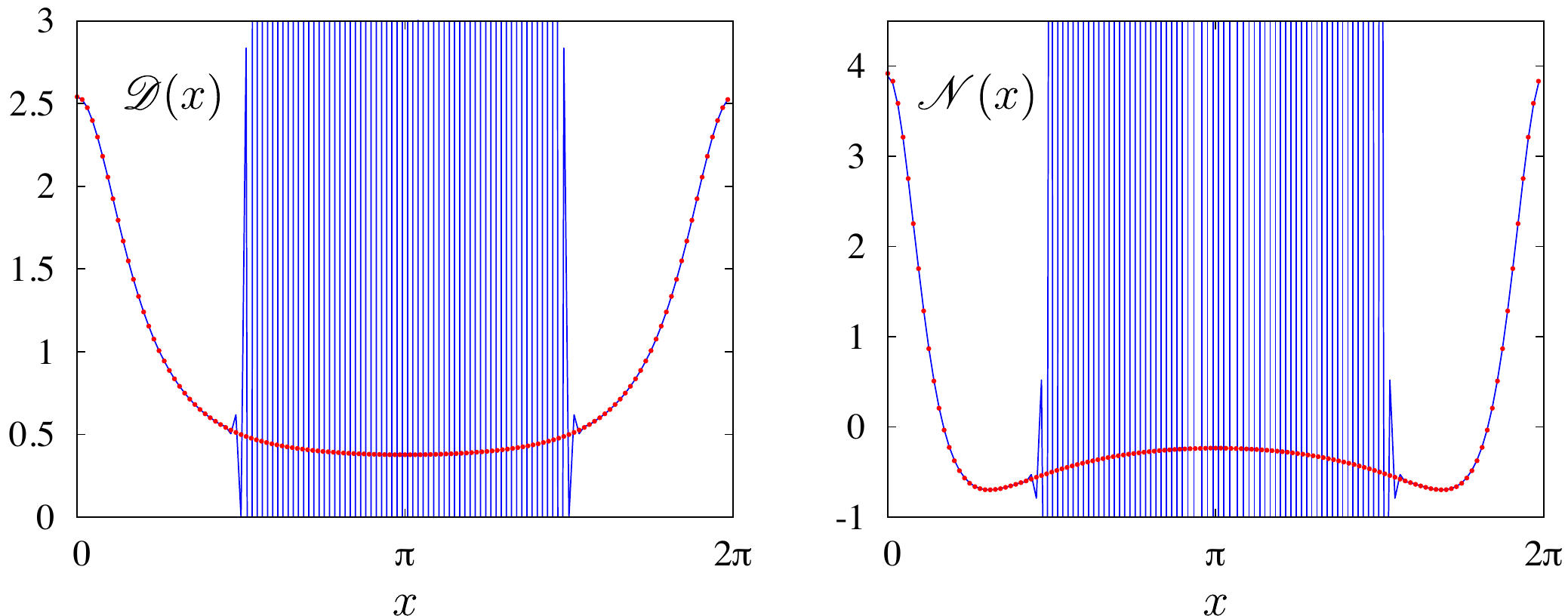}
\caption{\label{fig:wilken} Comparison of the series
  representation (\ref{Dirichlet:h:inf}), (\ref{Normal:h:inf}),
  (\ref{eq:ck:exact}) of
  $\sD$ and $\sN$, (blue curves), to the exact formulas
  (\ref{DN:exact}), (red markers), for $\eta(x)=-\veps\cos x$,
  $\veps=0.5$.  The series was truncated at $|k|=64$, and the results
  were plotted at 128 grid points.  A centered finite
  difference was used to compute the derivative in
  (\ref{Normal:h:inf}). As expected, the series diverges for
  $\frac{\pi}{2}\le x\le\frac{3\pi}{2}$, where $\eta(x)\ge0$.
}
\end{figure}

The simplest example illustrating these difficulties is the infinite
depth case with
\begin{equation}
  \eta(x)=-\veps\cos(x).
\end{equation}
Introducing the weights
$w_k=\cosh(k(\eta_\text{max}+h))$ in (\ref{wts}) and taking the limit as
$h\rightarrow\infty$, the system (\ref{AFM}) becomes
\begin{equation}\label{AFM:h:inf}
    \int_0^L e^{|k|(\eta-\eta_\text{max})}e^{-ikx}\sN(x)\,dx =
    \int_0^L -i\opn{sgn}(k)e^{|k|(\eta-\eta_\text{max})}e^{-ikx}\partial_x\sD(x)\,dx,
\end{equation}
where $\eta_\text{max}=\max_{0\le x\le2\pi}\eta(x) = \veps$.  Similarly, (\ref{Dirichlet}),
(\ref{Normal}) become
\begin{align}
  \label{Dirichlet:h:inf}
  &\sum_k c_k e^{|k|(\eta-\eta_\text{max})}e^{ikx} = \sD, \\
  \label{Normal:h:inf}
  &\sN = (-i\partial_x)\sum_k c_k\opn{sgn}(k)e^{|k|(\eta-\eta_\text{max})}e^{ikx}.
\end{align}
In our case, using $\cot\frac{x+iy}{2} = \frac{\sin x - i\sinh y}{\cosh y - \cos x}$
in (\ref{wilken_ex}), we have
\begin{align}\notag
  \phi(x,y) &= \frac{1}{2}\left(\frac{-\sinh y}{\cosh y-\cos x}+1\right), \quad
  \sD(x) =
  \frac{1}{2}\left(\frac{\sinh(\veps\cos x)}{\cosh(\veps\cos x)-\cos x}+1\right), \\
  \label{DN:exact}
  \sN(x) &= \phi_y - \eta_x\phi_x =
  \frac{\cosh(\veps\cos x)\cos x-1 +
    \veps (\sin^2x)\sinh(\veps\cos x)}{2(\cosh(\veps \cos x)-\cos x)^2}.
\end{align}
%
Since $\eta$ dips below the poles at $x\in2\pi\mathbb{Z}$, $\phi$ is
harmonic on $-\infty<y<\eta(x)$. Moreover, $\sD(x)=\phi(x,\eta(x))$ is
real analytic and $2\pi$-periodic.  Nevertheless, there is no solution
of (\ref{Dirichlet:h:inf}) valid over the whole interval $0\le x\le
2\pi$.  From (\ref{phi:ex1}), we see that the coefficients
\begin{equation}\label{eq:ck:exact}
  c_k = \begin{cases}
    1 & k=0 \\
    (1/2)e^{|k|\eta_\text{max}} & k\ne0
    \end{cases}
\end{equation}
will work over $\{x\,:\,\eta(x)<0\}$, but not elsewhere; see
Figure~\ref{fig:wilken}.

There remains the possibility that for any prescribed tolerance, a
choice of the $c_k$ can be made such that (\ref{Dirichlet:h:inf}) is
satisfied approximately, to the specified tolerance.  In other words,
the left-hand side is not treated as a series, but instead as a finite
linear combination of basis functions that can approximate $\sD$ to
arbitrary accuracy.  Rather than just add more terms to improve
accuracy, it may be necessary to change all the coefficients $c_k$.
To explore this possibility, we construct the $M$ by $K-1$ matrices
$A$ and $B$ with entries
\begin{equation}\label{eq:AB}
  A_{jk} = \frac{1}{M}\exp\Big(|k|\big[\eta(x_j)-\eta_\text{max}\big]\Big)
  e^{ikx_j}, \qquad
  B_{jk} = i\opn{sgn}(k)A_{jk},
\end{equation}
where $-K/2<k<K/2$, $x_j=2\pi j/M$, and $M\ge K$.  We then compute the
singular value decomposition $A=USV^*$, with $S$ and $V$ square and
$U$ of dimension $M\times(K-1)$, and evaluate
\begin{alignat}{2}
  \label{N:AFM}
  \sN &= U\opn{pinv}(S)V^*B^*\partial_x \sD, & &\text{(AFM)}, \\
  \label{N:AFM1}
  \sN &= -\partial_x B V\opn{pinv}(S)U^* \sD, & \quad\;\; &\text{(AFM$^*$)}.
\end{alignat}
Here $\partial_x$ is computed spectrally via the FFT (with no filter).
The idea here is to sample each column of $A$ and $B$ in (\ref{wts}) with
enough points that computing the SVD of the resulting matrix $A$ in
(\ref{eq:AB}) is equivalent (up to scaling by $\sqrt{M}$ in various
places) to computing the SVD of the quasi-matrix $A$ in (\ref{wts}),
whose columns are continuous, $L^2$ functions:
\begin{equation}\label{eq:quasi}
  A = USV^*, \qquad U:\mathbb{C}^{K-1}\rightarrow L^2, \qquad
  S:\mathbb{C}^n\rightarrow\mathbb{C}^n, \qquad
  V:\mathbb{C}^n\rightarrow\mathbb{C}^n.
\end{equation}
Here $A$ has been truncated to have $K-1$ columns, $U$ and $V$ are
unitary, $S$ is diagonal with positive decreasing entries, and $L^2$
is equipped with the inner product $\langle f,g\rangle=
\frac{1}{2\pi}\int_0^{2\pi} f\bar g\,dx$ to avoid factors of
$\sqrt{2\pi}$ elsewhere.  The columns of $U$ form an orthonormal basis
for the column-span of $A$.
Once sufficient grid resolution is reached, approximating $\sD$ and
$\sN$ at the $M$ collocation points leads to accurate approximation
throughout $(0,2\pi)$, using trigonometric polynomials to interpolate
between grid points. The columns of the matrix version of $U$ can be
thought of as sampled versions of the columns of the quasi-matrix $U$
from (\ref{eq:quasi}), up to a factor of
$\sqrt M$.
For smaller values of $M$, trigonometric interpolation becomes less
accurate, and the errors can be amplified significantly on division by
small singular values.  We remark that if the pseudo-inverse in
(\ref{N:AFM}) and (\ref{N:AFM1}) is replaced by an inverse, then
one can use a $QR$ factorization instead of the SVD to obtain
an orthonormal basis $U$ for the column span of $A$.  This gives
up some flexibility in regularizing the AFM/AFM$^*$ methods, but
is cheaper and has the advantage that the leading basis functions
do not change if $K$ is increased.  This $QR$ approach is equivalent
to Gram-Schmidt orthogonalization when the columns are sampled
sufficiently.

\begin{figure}
\includegraphics[width=\linewidth,trim=30 0 10 -10]{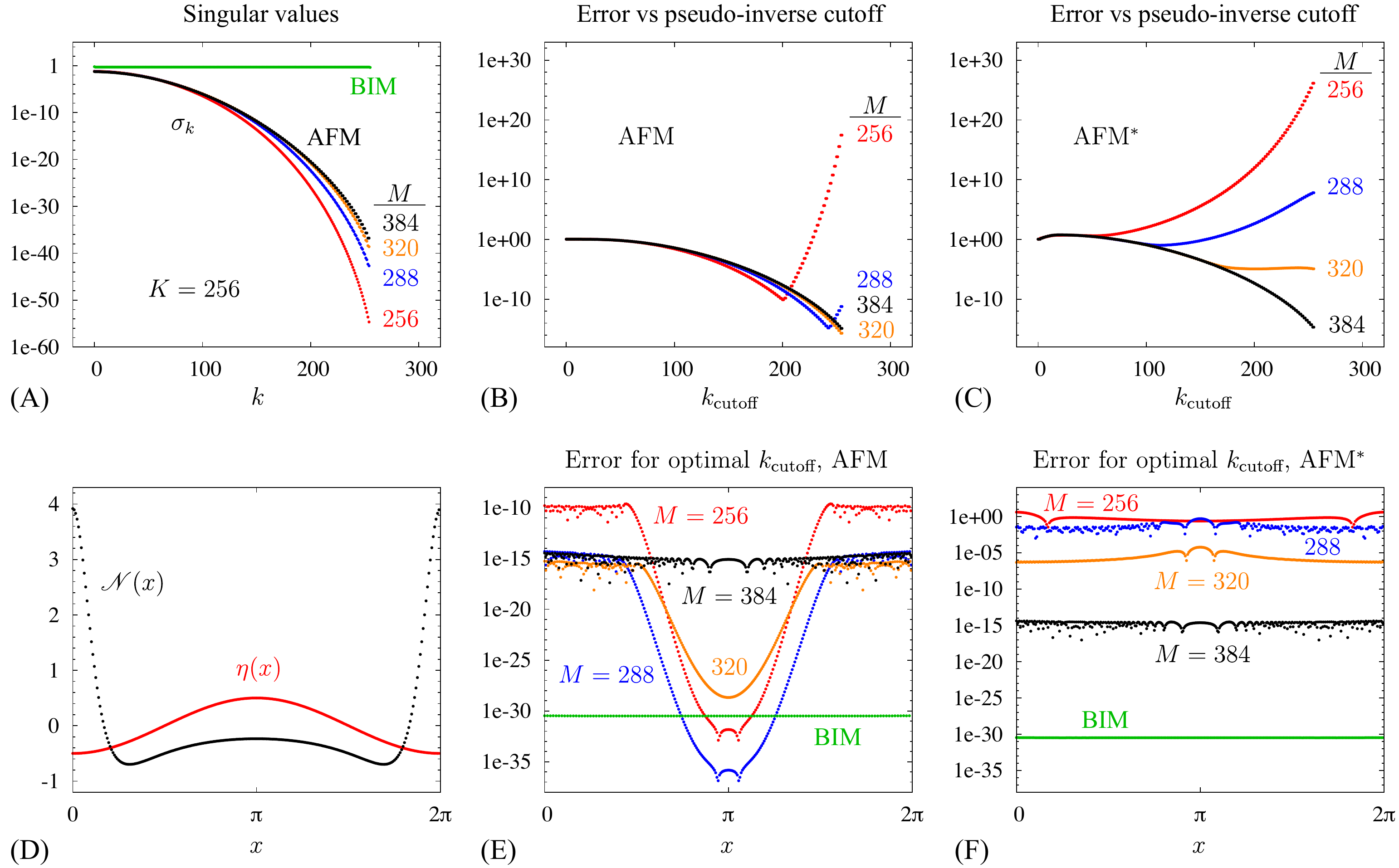}
\caption{\label{fig:off0} Singular values (A), exact solution (D), and
  effect of varying $M$ and the pseudo-inverse cutoff,
  $k_\text{cutoff}$, in the AFM and AFM$^*$ methods (B,C,E,F).  Errors
  were computed relative to the exact solution $\sN$ given in
  (\ref{DN:exact}). 360 bits of precision were used in the AFM and
  AFM$^*$ calculations so that aliasing and truncation errors dominate
  roundoff errors.  }
\end{figure}

The results of computing $\sN$ from $\sD$ in (\ref{DN:exact}) are
shown in Figure~\ref{fig:off0}.  The first panel shows the singular
values of $A$ with $K=256$ and $M\in\{256, 320, 288, 384\}$.  In all
four cases, the singular values decay very rapidly, with slightly
slower decay when $M$ is larger. By contrast, the singular values of
$A$ in the boundary integral approach (BIM) remain nearly constant.
Panel B shows the error in $\sN$ from (\ref{N:AFM}),
\begin{equation}\label{eq:error}
  \text{error} = \sqrt{\jt\frac{1}{M}\sum_{j=0}^{M-1} |E_j|^2}, \qquad
  E_j = \sN(x_j) -
  [U\opn{pinv}(S)V^*B^*\partial_x \sD]_j,
\end{equation}
where $\sN$ refers to the exact solution (\ref{DN:exact}), plotted in
Panel D. The error depends on the pseudo-inverse cutoff index,
$k_\text{cutoff}$, defined by
\begin{equation*}
  \opn{pinv}(S)_{ij} = \begin{cases} S_{ii}^{-1} & i=j\le k_\text{cutoff}, \\
    0 & \text{otherwise}.  \end{cases}
\end{equation*}
It consists of two parts, one due to how well $\sN$ is approximated by
the leading $k_\text{cutoff}$ columns of $U$, and one by how well the
coefficients $c=UU^*\sN$ are approximated by
\begin{equation*}
  c \approx \opn{pinv}(S)V^*B^*\partial_x \sD.
\end{equation*}
For all values of $M\ge K$, the error in panel B decreases initially
as $k_\text{cutoff}$ increases. This suggests that the leading entries
of $\opn{pinv}(S)V^*B^*\partial_x \sD$ are a good approximation of those
of $c$.  When $M$ is close to $K$ ($M=256$ or $288$ in the figure),
the error reaches a minimum at an optimal $k_\text{cutoff}$, and then
increases rapidly.  This occurs because the singular values $\sigma_k$
become so small that the high-index entries of
$\opn{pinv}(S)V^*B^*\partial_x \sD$ become large and no longer
approximate the corresponding entries of $c$.  For larger $M$, the
singular values decay more slowly and the error curve decreases
monotonically all the way to $k_\text{cutoff}=K-1$.

Panel E shows the vector version of the error, namely $E_j$ in
(\ref{eq:error}), corresponding to the minima of the error curves in
Panel B. For smaller values of $M$, the error $E_j$ is largest where
$\eta(x)$ is smallest. This is not surprising since the corresponding
rows of $A$ are smaller in this region due to the exponential growth
of the basis functions $e^{ky}e^{ikx}$ in the $y$-direction. What is
surprising is the extreme accuracy that is achieved by the AFM method
in the region where $\eta(x)>0$ for smaller values of $M$.  We do not
know why this occurs. This additional accuracy disappears as the
errors are reduced in the region where $\eta(x)<0$ by increasing $M$.
Once $M$ reaches 384, the columns of $A$ are well-resolved in
$L^2(0,2\pi)$ as discussed above, and the error $E_j$ is roughly
uniform throughout the domain.

Panels C and F show the same results for the AFM$^*$ method.  The
results are very poor when $M=K=256$, presumably a consequence of
aliasing errors in sampling the columns of $A$ being amplified on
division by small singular values.  As $M$ increases, the aliasing
errors become smaller and the singular values become larger.  By the
time $M$ reaches 384, the error of the AFM$^*$ method is similar to
that of the AFM method, of order $10^{-15}$.  By contrast, the
boundary integral method has errors of order $10^{-30}$ with
$K=M=256$.  There is little amplification of roundoff error in the BIM
approach since $A$ is so well-conditioned.

\begin{figure}[t]
\includegraphics[width=\linewidth,trim=30 0 10 -10]{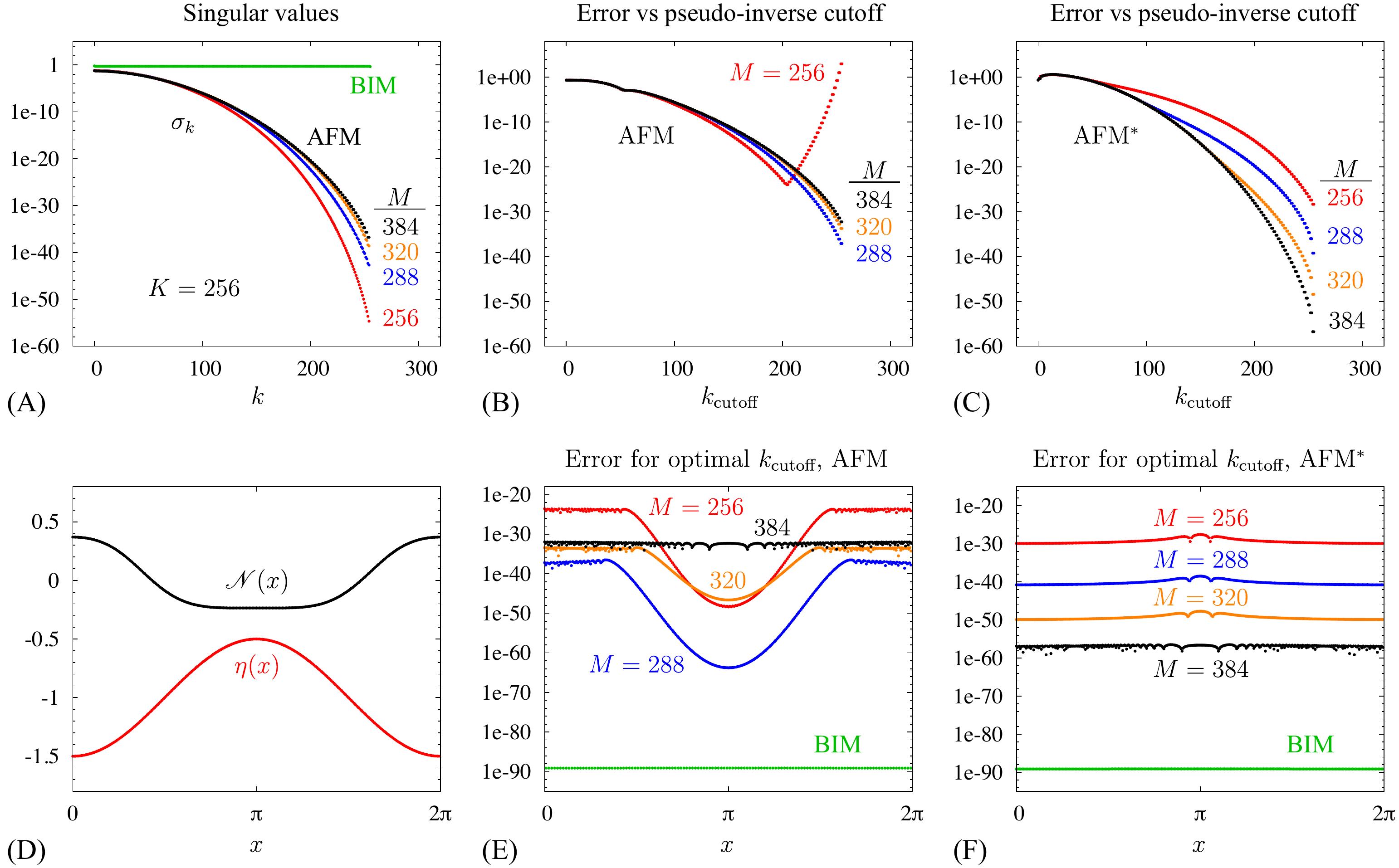}
\caption{\label{fig:off1} Same as Figure~\ref{fig:off0}, but with a
  modified wave profile $\eta(x)$.  }
\end{figure}

Figure~\ref{fig:off1} shows a similar computation to the above,
but with $\eta(x)$ offset vertically so that the series
(\ref{phi:ex1}) converges at all points on the curve.
Specifically, we set
\begin{equation*}
  \eta(x) = -1.0 - \veps\cos(x), \qquad
  \veps = 0.5.
\end{equation*}
The formulas (\ref{DN:exact}) remain nearly the same, with $\veps\cos
x$ replaced by $1+\veps\cos x$.  The singular values in Panel A are
the same as those in Figure~\ref{fig:off0} since $\eta_\text{max}$ is
also shifted downward by 1 in (\ref{eq:AB}).  The AFM method behaves
similarly to before, achieving exceptional accuracy in the region
where $\eta$ is largest when $M$ is close to $K$, and achieving nearly
uniform accuracy once $M$ is large enough to fully resolve the columns of
the continuous version of $A$.  Because $\sD$ and $\sN$ are smoother
(with faster decay of Fourier modes), all three methods (AFM, AFM$^*$
and BIM) yield smaller errors in Figure~\ref{fig:off1} than in
Figure~\ref{fig:off0}. The AFM$^*$ method turns out to be superior to
the AFM method on this example for all four choices of $M$.  By
contrast, in Figure~\ref{fig:off0}, the AFM method was better for
$M=256, 288$ and $320$, and the two methods were equal when $M=384$.
A partial explanation is that the left-hand side of
(\ref{Dirichlet:h:inf}) is a genuinely convergent series in this
second example, as opposed to a finite linear combination of vectors
from a dense set.  Thus, the AFM$^*$ method has an easier time
selecting the coefficients $c_k$ that best represent $\sD$.  By
contrast, for the AFM method, there is little difference between the
two examples.  In either case, the exact solutions $\sD$ and
$\sN$ satisfy the AFM global relation, and the question of convergence
comes down to how well the trapezoidal rule approximates the integrals,
and how much the errors are amplified by the poorly-conditioned $A$
matrix.  We do not know why the AFM$^*$ method turns out to be
25 orders of magnitude more accurate once $M$ reaches 384. As before,
the BIM method is superior to both AFM methods.

\begin{figure}[t]
\includegraphics[width=\linewidth,trim=30 0 10 -10]{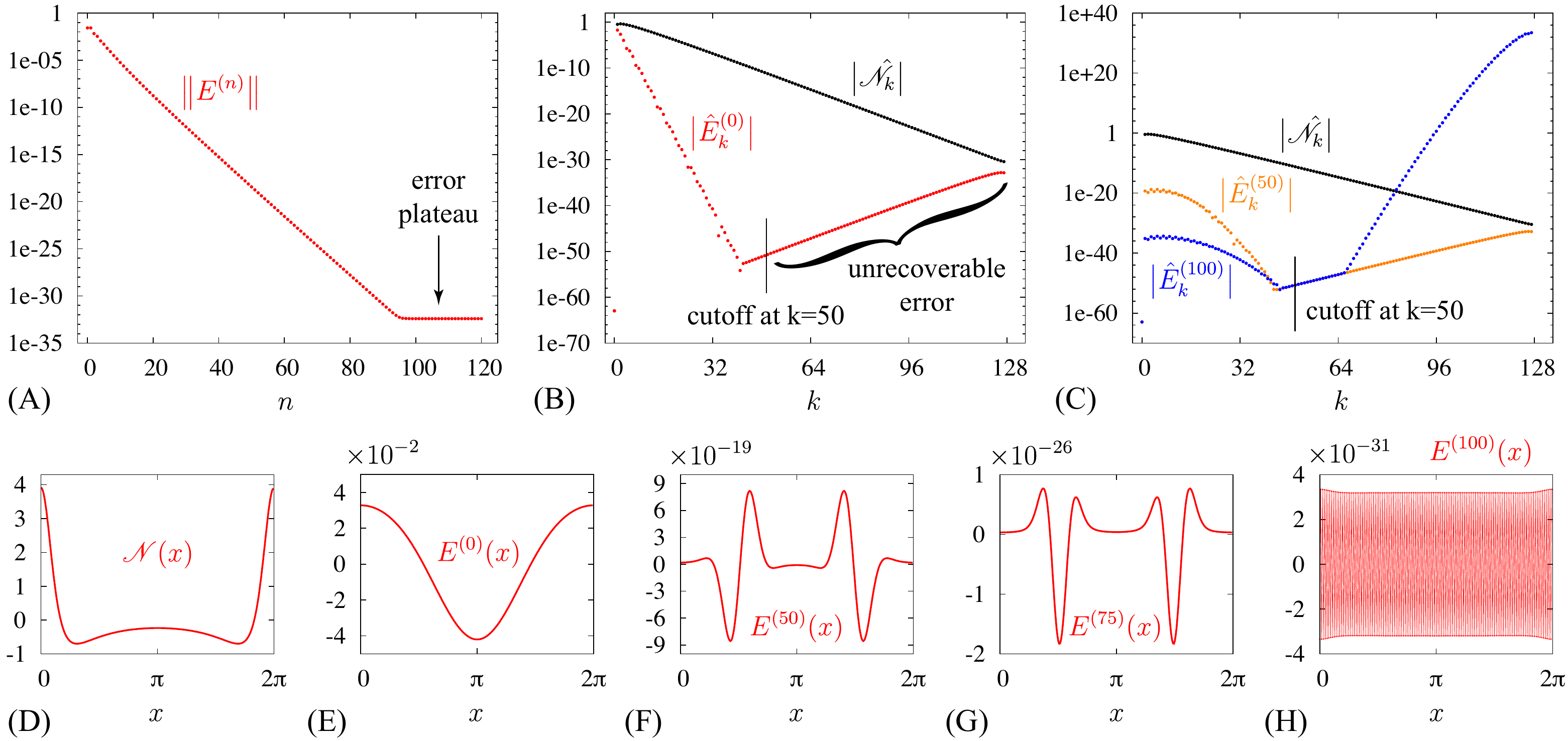}
\caption{\label{fig:craigErr} Error in the Craig-Sulem expansion at
  various orders. We omit the negative index Fourier modes since
  $\hat{\sN}_{-k}=\overline{\hat\sN_k}$ and $\hat E^{(n)}_{-k} =
  \overline{\hat E^{(n)}_k}$.  }
\end{figure}

In Figure~\ref{fig:craigErr}, we check convergence of the Craig-Sulem
expansion for this example.  As in Figure~\ref{fig:off0}, we consider
$\eta(x)=-\veps\cos x$, $\veps=0.5$, with Dirichlet data as in
(\ref{DN:exact}). The computations were done on a 256-point grid with
360 bits of precision. Panel A shows the errors in the partial
sums, defined as
\begin{equation}\label{eq:En:def}
  \big\|E^{(n)}\|=\sqrt{\jt\frac{1}{M}\sum_{j=0}^{M-1}
    \big|E^{(n)}(x_j)\big|^2},
  \qquad E^{(n)}(x) = \sN(x) - \sum_{j=0}^n [G_j(\eta)\sD](x).
\end{equation}
The errors decrease steadily until $n=95$, where $\big\|E^{(n)}\big\|$
reaches $10^{-32}$, the level of aliasing errors in the Fourier
representation of $\sD$ and $\sN$ on a 256-point grid.  Recall that
the error in the boundary integral method with 256 points was also
around $10^{-32}$.  Panel B shows the Fourier spectrum of $\sN(x)$
(from the exact solution) and $E^{(0)}(x)$.  The rapid decay of the
Fourier modes of $E^{(0)}(x)$ show that the zeroth order approximation
does an excellent job of predicting the high-frequency components of
$\sN(x)$, but not the low-frequency ones.  All further corrections
will only be made to the first 50 Fourier modes.  In other words, we
actually use
\begin{equation}
  \hat E_k^{(n)} = \begin{cases}
     \hat\sN_k - \sum_{j=0}^n [G_j(\eta)\sD]^\wedge_k, & |k|<50, \\
     \hat\sN_k - [G_0(\eta)\sD]^\wedge_k, & |k|\ge50
   \end{cases}
\end{equation}
when reconstructing $E(x)$ in Panels E--H.  As a result, the modes
marked ``unrecoverable error'' in Panel B are frozen, and will not be
altered by successive corrections.  This is reasonable since these
modes are dominated by aliasing errors in sampling $\sN(x)$ on the
256-point grid.  On a larger grid, if more accuracy were desired, the
cutoff would need to be increased. The error plateau in Panel A is due
to these aliasing errors.

Panel C shows the magnitude of $\hat E^{(n)}_k$ versus $k$ for $n=50$
and $n=100$.  Each successive correction reduces the error in the
leading Fourier modes of the Neumann data.  Note that the bulge near
$k=0$ in $\hat E^{(100)}_k$ is below $10^{-32}$, and therefore
additional corrections do not improve the global error, which is
dominated by high-frequency modes at that point.  The rapid growth of
$\hat E^{(100)}_k$ above $k=65$ is due to truncating the Fourier
series to $|k|\le128$ when computing the Craig-Sulem expansion.
Analogous behavior was seen in Figures~\ref{fig:cosB}--\ref{fig:sym},
where we observed that errors in the columns of $G_n(f)^\wedge$
propagate inward from high to low wave numbers.  Due to the cutoff,
these large errors do not affect the reconstruction of $\sN(x)$.
Panel D shows the exact solution $\sN(x)$ while panels E--H show the
errors in various partial sum reconstructions.  Because the low-order
modes are the least accurate (as seen in Panel C), the corrections are
remarkably smooth, non-oscillatory curves.  The highly oscillatory
error in $E^{(100)}$ is due to the ``unrecoverable error'' in Panel B;
if $M$ were increased, $E^{(100)}$ would also be very smooth.  In
summary, the Craig-Sulem expansion performs well on this example,
converging to the exact solution, up to aliasing errors of order
$10^{-32}$, in 95 iterations.

\subsection{Water wave examples.}\label{waterwave}
We conclude with two examples in which the wave profile and velocity
potential come from solutions of the water wave problem.  The first is
a large-amplitude standing water wave shortly before reaching maximum
height.  The second consists of two Stokes waves of different
amplitudes traveling to the right.  Both examples can be evolved
efficiently using the boundary integral method \cite{trav:stand,
  water2}.  Our interest here is whether the DNO expansion method and
the two AFM methods can take data ($\eta$ and $\sD$) that are only
known to double- or quadruple-precision accuracy and return Neumann
data ($\sN$) with the same accuracy.  We allow ourselves to use
additional precision for intermediate calculations.

\begin{figure}
\includegraphics[width=\linewidth,trim=30 0 10 -10]{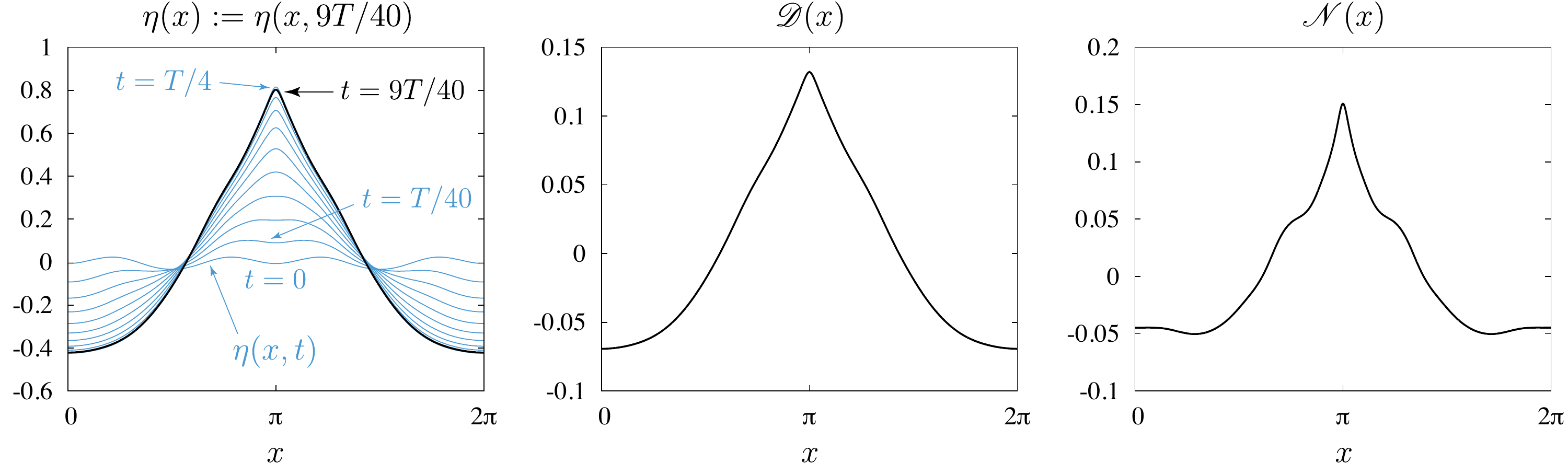}
\caption{\label{fig:stand2} Evolution of $\eta(x,t)$ over a quarter
  period, $T/4$, in increments of $T/40$, along with the Dirichlet and
  Neumann data corresponding to $t=9T/40$.  }
\end{figure}

\begin{figure}
\includegraphics[width=\linewidth,trim=30 0 10 -10]{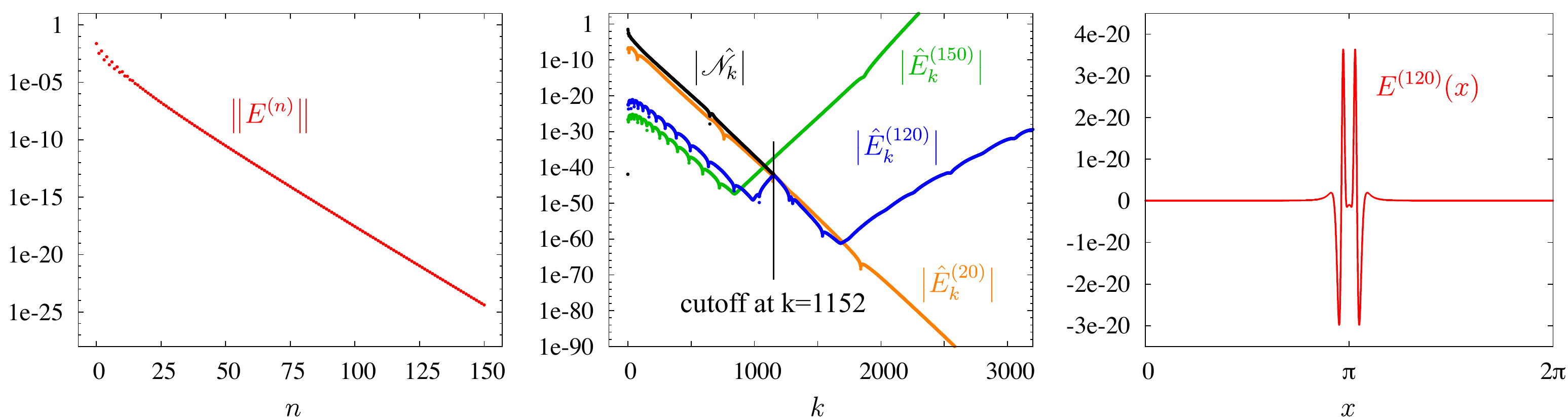}
\caption{\label{fig:craig:stand}
  Errors in the Craig-Sulem calculation of the Neumann data
  for the standing wave at $t=9T/40$. Fourier modes above $k=1152$
  are set to zero in the reconstruction of $\sN(x)$ and the computation
  of $E^{(n)}(x)$.
}
\end{figure}

\begin{figure}
\includegraphics[width=\linewidth,trim=30 0 10 -10]{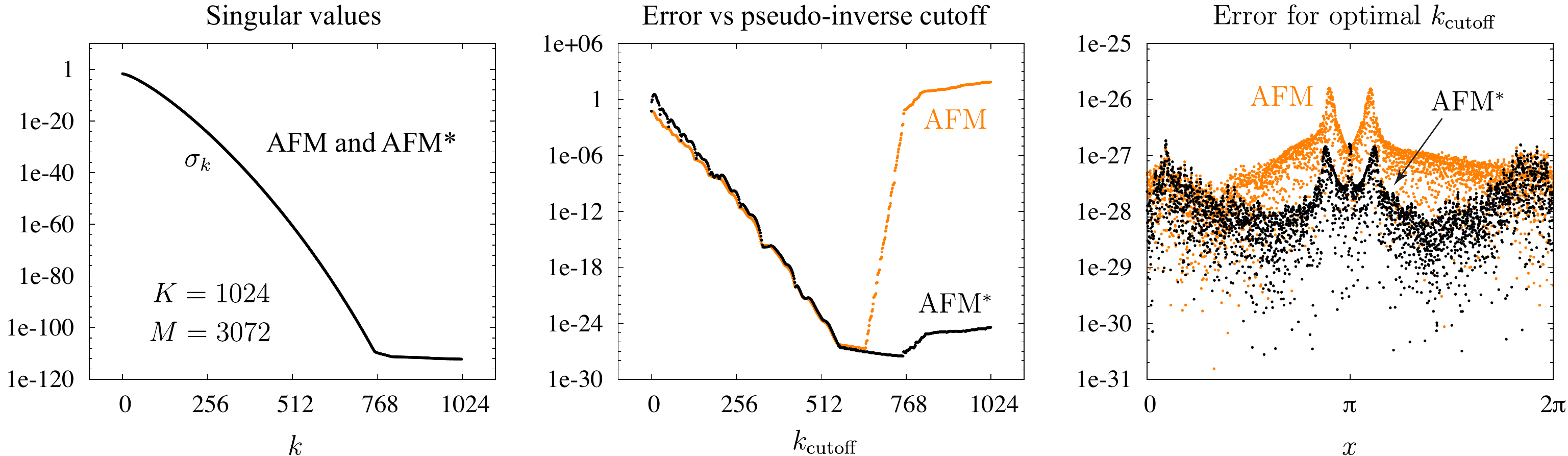}
\caption{\label{fig:afm:stand}
  AFM and AFM$^*$ calculation of the Neumann data for the
  standing wave at $t=9T/40$. The optimal cutoffs in the right panel
  correspond to the error minima in the center panel:
  $k_\text{cutoff}=643$ (AFM), $k_\text{cutoff}=759$ (AFM$^*$).
}
\end{figure}

Figures~\ref{fig:stand2}--\ref{fig:afm:stand} show the calculation of
the DNO operator for a large-amplitude standing water wave.  We
selected the infinite depth wave corresponding to the first local
maximum of wave height (half the vertical crest-to-trough distance)
for this example.  When wavelength is set to $L=2\pi$ and the
acceleration of gravity is $g=1$, this wave has period $T=6.53996$,
wave height $h=0.620173$, crest acceleration $A_c=0.926312$, and fifth
Fourier mode of $\varphi$ at $t=0$ (a good bifurcation parameter in
this regime) of $\hat\varphi_5(0) = 0.00245499$.  See
\cite{water1,water2} for details on how the wave was
computed. Since the wave comes to rest at $t=T/4$, $\sD\equiv0$ at
that time.  Thus, to avoid a trivial DNO calculation while keeping
$\eta$ close to its maximum-amplitude state, we selected the wave
profile at $t=9T/40$ for this example.  Figure~\ref{fig:stand2} shows
snapshots of $\eta(x,t)$ at equal intervals of size $\Delta t=T/40$
over a quarter-period, along with the Dirichlet and Neumann data
corresponding to $t=9T/40$.

Figure~\ref{fig:craig:stand} shows the results of the DNO calculation
using the Craig-Sulem method. The left panel shows the error as a
function of order for $0\le n\le 150$.  The middle panel shows the
Fourier decomposition of the error at orders $n=20$, $120$ and $150$,
together with the Fourier modes of the ``exact'' solution.  Unlike the
results of Figure~\ref{fig:craigErr}, the zeroth order approximation
does not lead to a large improvement in the high-frequency modes (seen
as a faster decay rate in $|\hat E_k^\e0|$ than $|\hat\sN_k|$ in
Figure~\ref{fig:craigErr}).  In fact, $|\hat E_k^\e 0|$ is difficult
to distinguish from $|\hat \sN_k|$, so we plotted $|\hat E_k^\e{20}|$
instead, which has a similar decay rate but is shifted down slightly.
The data for this problem, $\eta(x)$ and $\sD(x)$, are specified via
their leading 750 Fourier modes (recorded in quadruple-precision,
i.e.~32 digits).  This gives an approximation of the standing wave to
25 digits of accuracy, but is regarded here as specifying the DNO
problem with infinite precision.  The ``exact'' solution was computed
using the boundary integral method with 2304 collocation points and
212 bits of precision, leading to approximately 40 digits of accuracy.
The resulting 1152 Fourier modes of this ``exact'' solution are
labeled $\hat\sN_k$ in the middle panel. The Craig-Sulem expansion was
performed using $M=16384$ grid points for the FFT and 900 bits of
precision in intermediate calculations in order to achieve accurate
results up to order $n=150$ for modes $|k|<K/2$, $K=2304$. This cutoff
$|k|<1152$ was chosen to agree with the last computed mode of the
``exact'' solution.  The right panel shows the error in the 120th
order approximation as a function of $x$. We see that the error is a
smooth function of order $10^{-20}$ that is largest near the crest tip
at $x=\pi$.  In summary, using 900 bits (271 digits) and $M=16384$
grid points, the 150th order Craig-Sulem method is able to achieve
errors around $10^{-25}$, which is comparable to the original standing
wave calculation, which was done in quadruple-precision (32 digits)
with 2048 grid points using the boundary integral method.  Though the
Craig-Sulem method is not competitive, it is still remarkable that
such a large-amplitude wave would be inside the radius of convergence
of the DNO expansion.

Figure~\ref{fig:afm:stand} shows the same calculation using the AFM
and AFM$^*$ methods.  The same ``exact'' solution as in
Figure~\ref{fig:craig:stand} was used to measure errors.  The large
(crest-to-trough) wave height leads to very rapidly-decaying singular
values.  Nevertheless, using $K=1024$ modes and $M=3072$ collocation
points with 360 bits (108 digits) of precision, we are able to achieve
25 digits of accuracy in the solution.  The AFM and AFM$^*$ methods
are comparable for cutoffs up to about $k_\text{cutoff}=640$. After
that, the error in the AFM method grows rapidly while the AFM$^*$
method flattens out. Thus, the AFM$^*$ method is somewhat more robust
in this example.

\begin{figure}
\includegraphics[width=\linewidth,trim=30 0 10 -10]{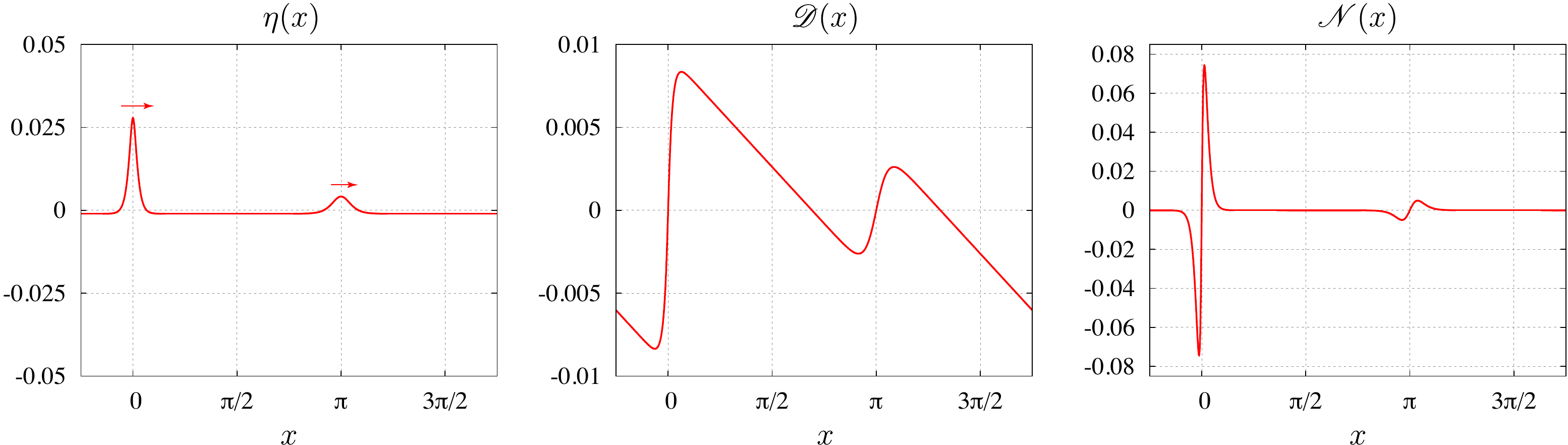}
\caption{\label{fig:stokes2}
  Plots of $\eta(x)$, $\sD(x)$ and $\sN(x)$ for the
  superposition of Stokes waves with parameters in (\ref{eq:stokes2})
  and fluid depth $h=0.05$.
}
\end{figure}

\begin{figure}
\includegraphics[width=\linewidth,trim=30 0 10 -10]{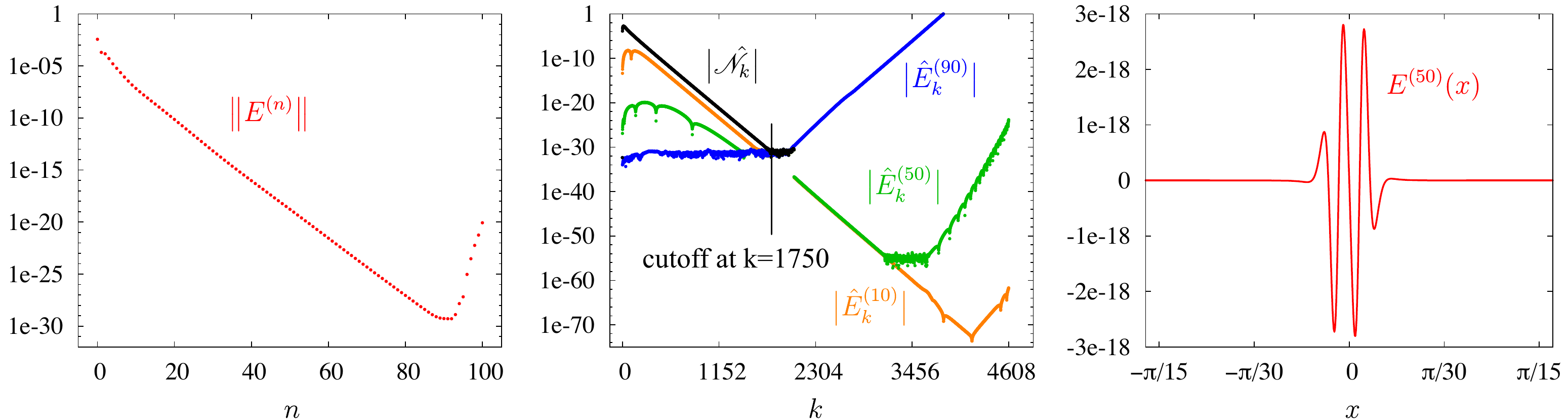}
\caption{\label{fig:craigStokes}
  Errors in the Craig-Sulem calculation of the Neumann data for a
  pair of traveling Stokes waves with parameters (\ref{eq:stokes2})
  and fluid depth $h=0.05$.
}
\end{figure}

The final example of this section consists of two Stokes waves
traveling right on a $2\pi$-periodic domain.  Initially, one wave is
centered at $x=0$ and the other is centered at $x=\pi$.  The fluid
depth is $h=0.05$ and the parameters of the waves are
\begin{equation}\label{eq:stokes2}
  \begin{array}{r|ccc}
    & \hat\eta_1 & \eta_\text{max}-\eta_{min} & c \\ \hline
    \text{wave 1} & \phantom{-}7.4\times 10^{-4} & 0.028919 & 0.27349 \\
    \text{wave 2} & -3.0\times 10^{-4} & 0.005202 & 0.23290
  \end{array}
\end{equation}
where $c$ is the wave speed.  Plots of $\eta(x)$, $\sD(x)$ and
$\sN(x)$ are given in Figure~\ref{fig:stokes2}.
Figure~\ref{fig:craigStokes} shows the results of the DNO calculation
using the Craig-Sulem method.  The recursion (\ref{recur}) is
modified as follows to account for finite depth:
\begin{align*}
  G_0(f) &= |D|\tanh(h|D|), \qquad
  Y_s = \left\{\begin{array}{ll}
      |D|^s, & s\text{ even} \\
      |D|^s\tanh(h|D|), & s\text{ odd}
    \end{array}\right\}, \\
  G_n(f) &= A_n(f) - \sum_{s=1}^{n-1}\frac{1}{(n-s)!}Y_{n-s}f^{n-s}G_s(f), \qquad
  n=1,2,3,\dots, \\
  A_n(f) &= \left\{\begin{array}{ll}
      (n!)^{-1}|D|^{n-1}\big( \tanh(h|D|) D f^n D - |D| f^n G_0 \big) ,&
      n\text{ even} \\[2pt]
      (n!)^{-1}|D|^{n-1}\big( D f^n D - G_0 f^n G_0 \big), & n\text{ odd}
    \end{array}\right\}, \\
  A_n(f)^\wedge_{kj} &= \frac{jk^n}{n!} a_{nkj} (f^n)^\wedge_{k-j}, \qquad
  a_{nkj} = \left\{\begin{array}{ll}
      \tanh kh - \tanh jh, & n\text{ even} \\
      1 - (\tanh kh)(\tanh jh), & n\text{ odd}
      \end{array}\right\}.
\end{align*}
We used $M=9216$ grid points for the
FFT and computed the expansion through order $n=100$.  The ``exact''
solution was computed via the boundary-integral method in
quadruple-precision, which is correct to about 28 digits of
accuracy.  The Craig-Sulem expansion reaches this level of accuracy at
order $n=90$.  The errors then increase (left panel) due to the error
growth region crossing the cutoff mode $k=1750$ at 92nd order (center
panel).  Increasing $M$ would delay this crossing, but is unnecessary
since $90$th order is sufficient to reach the accuracy of the
underlying ``exact'' solution.  The right panel shows the error as a
function of $x$ at order $n=50$.  The error is concentrated near
$x=0$, where $\eta(x)$ is largest. In the region not shown, the error
is uniformly less than $10^{-29}$, even near $x=\pi$, the location of
the second Stokes wave. Thus, most of the work goes into resolving the
solution near the larger peak.

\begin{figure}
\includegraphics[width=\linewidth,trim=30 0 10 -10]{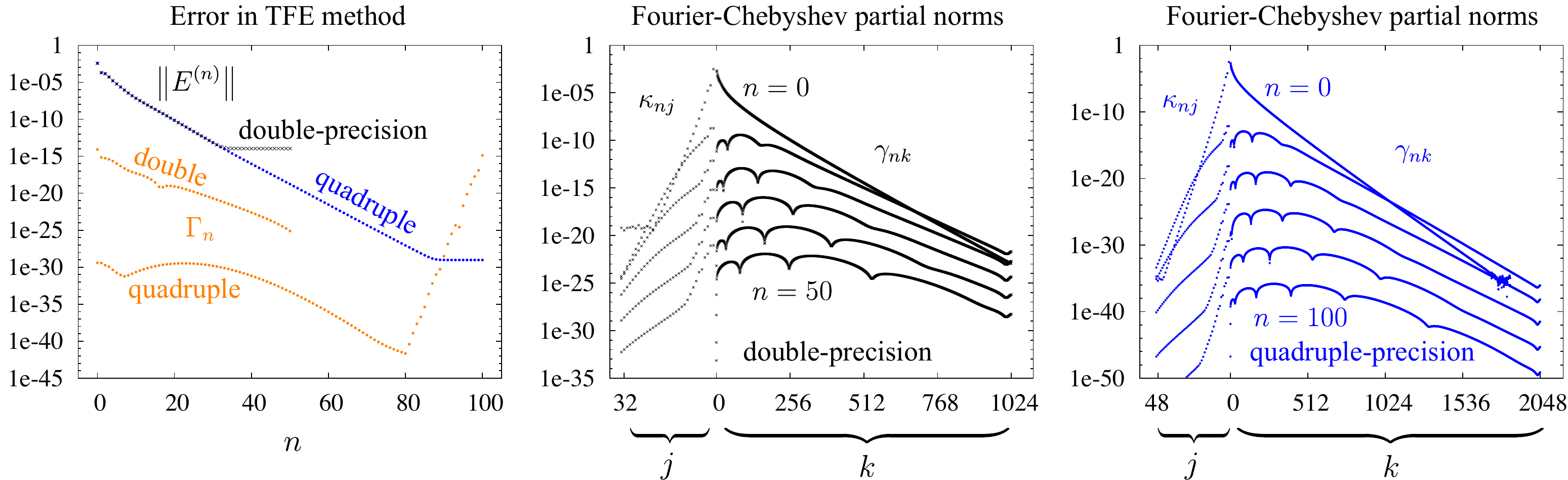}
\caption{\label{fig:tfeStokes} Two measures of error ($\|
  E^{(n)}\|$ and $\Gamma_n$) and Fourier-Chebyshev partial norms
  (\ref{eq:four:cheb:nrm}) for the TFE calculation of the Neumann data
  for a pair of traveling Stokes waves with parameters
  (\ref{eq:stokes2}) and fluid depth $h=0.05$.  }
\end{figure}

Since this last example has finite-depth, we can also use the variant
of the transformed field expansion method described in \S\ref{sec:tfe}
to compute the Neumann data from the Dirichlet data. The left panel of
Figure~\ref{fig:tfeStokes} shows the error $\|E^{(n)}\|$ from
(\ref{eq:En:def}) in the Neumann data as a function of the order $n$,
as well as the difference between the $n$th order terms from the
TFE and CS expansions,
\begin{equation}
  \Gamma_n = \big\| G_n^{TFE}(\eta)\sD - G_n^{CS}(\eta)\sD \big\|.
\end{equation}
The error $\|E^{(n)}\|$ reaches a plateau of $10^{-14}$ for $n\ge33$
in double-precision, and $10^{-29}$ for $n\ge88$ in
quadruple-precision.  As before, the ``exact'' solution was computed
in quadruple-precision using the boundary integral method. The orange
markers show that the terms in the TFE expansion agree with the
corresponding terms in the CS expansion to roundoff error accuracy. In
particular, a plot of $E^{(50)}(x)$ for the TFE method (not shown)
looks identical to that of the CS expansion in the right panel of
Figure~\ref{fig:craigStokes}.

We also note in Figure~\ref{fig:tfeStokes} that $\Gamma_n$ exhibits a
downward trend as $n$ increases, indicating that the terms in the TFE
expansion maintain several correct digits of relative accuracy beyond
the point that $\|E^{(n)}\|$ reaches the plateau region, i.e.~the
point where successive terms are smaller in magnitude than the
absolute errors of the leading terms.  In the double-precision case,
most of the error in the plateau region is due to the error in the
zeroth order term --- $\Gamma_0$ is larger than the sum of the other
$\Gamma_n$.  In quadruple-precision, $\Gamma_0$ is still largest,
although the $\Gamma_n$ with $15\le n\le 30$ are of comparable
size. This may be partly due to the cutoff at $k=1750$ used in the CS
expansion in Figure~\ref{fig:craigStokes} to eliminate high-frequency
noise. Certainly, the rapid growth in $\Gamma_n$ for $n\ge80$ in
Figure~\ref{fig:tfeStokes} is due to errors in the CS expansion rather
than the TFE expansion. Indeed, $\|E^{(n)}\|$ in
Figure~\ref{fig:craigStokes} grows rapidly for $n\ge90$ while it
remains flat for $n\ge90$ in Figure~\ref{fig:tfeStokes}.  Recall that
a third method (the BIM method) was used for the exact solution in both
plots.

The middle and right panels of Figure~\ref{fig:tfeStokes} were used
to decide how many grid points to use in the TFE method. We used
$M=2048$, $N=32$ in double-precision and $M=4096$, $N=48$ in
quadruple-precision.  These plots show the partial norms
\begin{equation}\label{eq:four:cheb:nrm}
  \kappa_{nj} = \sqrt{\jt \sum_k |\alpha^n_j(k)|^2}, \qquad
  \gamma_{nk} = \sqrt{\jt \sum_j |\alpha^n_j(k)|^2},
\end{equation}
where $\alpha^n_j(k)$ are the Chebyshev coefficients of the
$k$th Fourier mode of $u_n$:
\begin{equation*}
  \hat u_n(k,y) = \sum_{j=0}^N \alpha^n_j(k)\,T_j(1+2h^{-1}y), \qquad
  (|k|\le M/2,\; -h<y<0).
\end{equation*}
Recall that $u(x,y) = \phi(x,(1+h^{-1}\eta)y+\eta)$ is defined on a
rectangle and expanded in powers of $\veps$, $u(x,y) = \sum_n
\veps^nu_n(x,y)$, where $\eta=\veps f$. The partial norms
$\kappa_{nj}$ and $\gamma_{nk}$ are the norms of the rows and columns
of the matrix of Fourier-Chebyshev coefficients of the function
$u_n(x,y)$.  To resolve the solution spectrally, the mesh needs to be
large enough that the partial norms decay to the desired tolerance as
$k\rightarrow M/2$ or $j\rightarrow N$.  As shown in
Figure~\ref{fig:tfeStokes}, the meshes we selected are sufficient
to reach roundoff level tolerances in these limits.

It is worth noting that many fewer grid points are needed in the
$y$-direction than in the $x$-direction, so the price of discretizing
the bulk fluid is not as severe as one might imagine.  Moreover, as
with the BIM method, intermediate calculations can be done in double
or quadruple-precision arithmetic to achieve similar levels of
accuracy in the solution.  The overall running times (in seconds) of
the various methods on this example are given in the following table:
\begin{equation*}
  \begin{array}{r|c|c|c|c|c|c|c}
    \text{method} & \text{BIM(d)} & \text{BIM(q)} &
    \text{TFE(d)} & \text{TFE(q)} & \text{CS} &
    \text{AFM} & \text{AFM-QR} \\ \hline
    \text{time} & 0.156 & 3.53 & 2.46 & 24.0 & 2762 & 11222 &
    920
  \end{array}
\end{equation*}
The code was run on a 3.33 GHz Intel Xeon X5680 system with 12 cores.
Here (d) and (q) stand for double and quadruple-precision, and the
other methods were run with 360 bits (108 digits) of precision.  The
running time of AFM$^*$ is nearly identical to that of AFM, and the
AFM-QR variant will be described below.  Clearly, the BIM and TFE
approaches are superior to the CS and AFM-based methods since
arbitrary-precision arithmetic is not required.

\begin{figure}
\includegraphics[width=\linewidth,trim=30 0 10 -10]{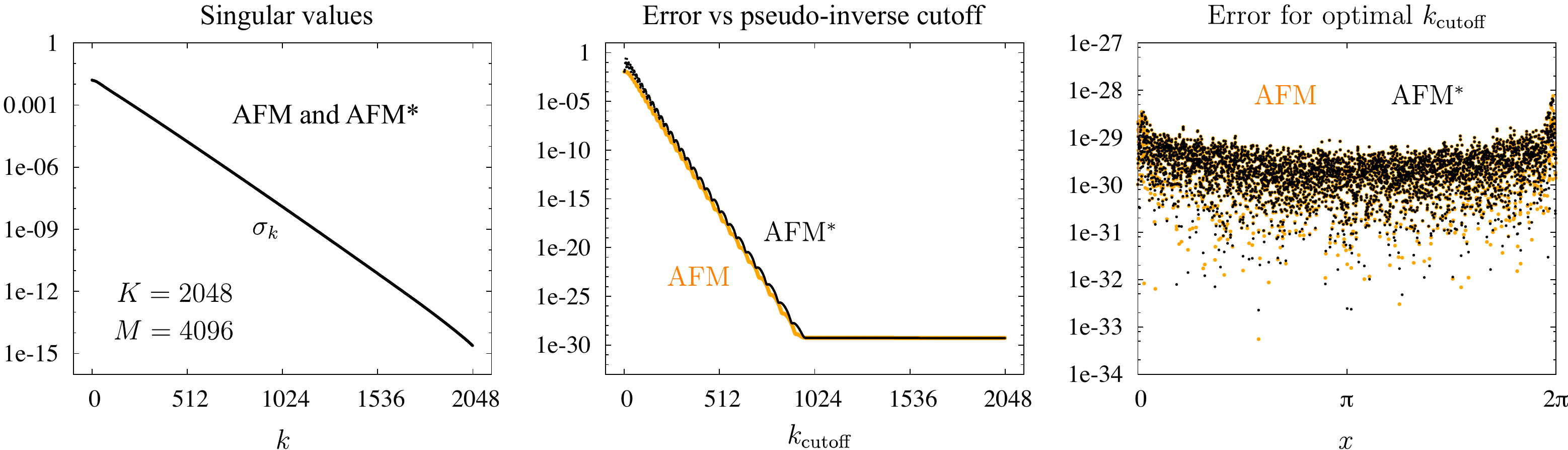}
\caption{\label{fig:afmStokes}
  AFM and AFM$^*$ calculation of the Neumann data for the
  superposition of Stokes waves in Figure~\ref{fig:stokes2}.  These
  AFM methods are much better conditioned in this shallow water regime
  than the deep water cases considered above.  The optimal cutoff
  modes at right were $k_\text{cutoff}=2046$ (AFM) and $2045$ (AFM$^*$).
}
\end{figure}

We conclude with the results of the AFM and AFM$^*$ methods on this
example, which we selected initially as being relevant to water
waves and likely to cause difficulties for the AFM and AFM$^*$
methods. Our reasoning was that the second wave is large enough to
require many Fourier modes to resolve its shape, but small enough that
$\cosh(k\eta(x))$ is many orders of magnitude larger at the crest of
the first wave than at the crest of the second (for large $k$).  Thus,
substantial cancellation must occur near $x=0$ in order to resolve the
behavior of $\phi$ near $x=\pi$.  However, the singular values in the
left panel of Figure~\ref{fig:afmStokes} decay slowly in comparison to
the infinite-depth standing wave case of Figure~\ref{fig:afm:stand}.
This is because the waves in these two examples have similar Fourier
decay rates for $\hat\eta_k$ and $\hat\sD_k$, but the deep-water
standing wave has a much larger vertical crest-to-trough height.  To
the extent that this is generally the case for waves in deep versus
shallow water, the AFM and AFM$^*$ methods appear to be better suited
for finite depth problems.

\begin{figure}
\includegraphics[width=\linewidth]{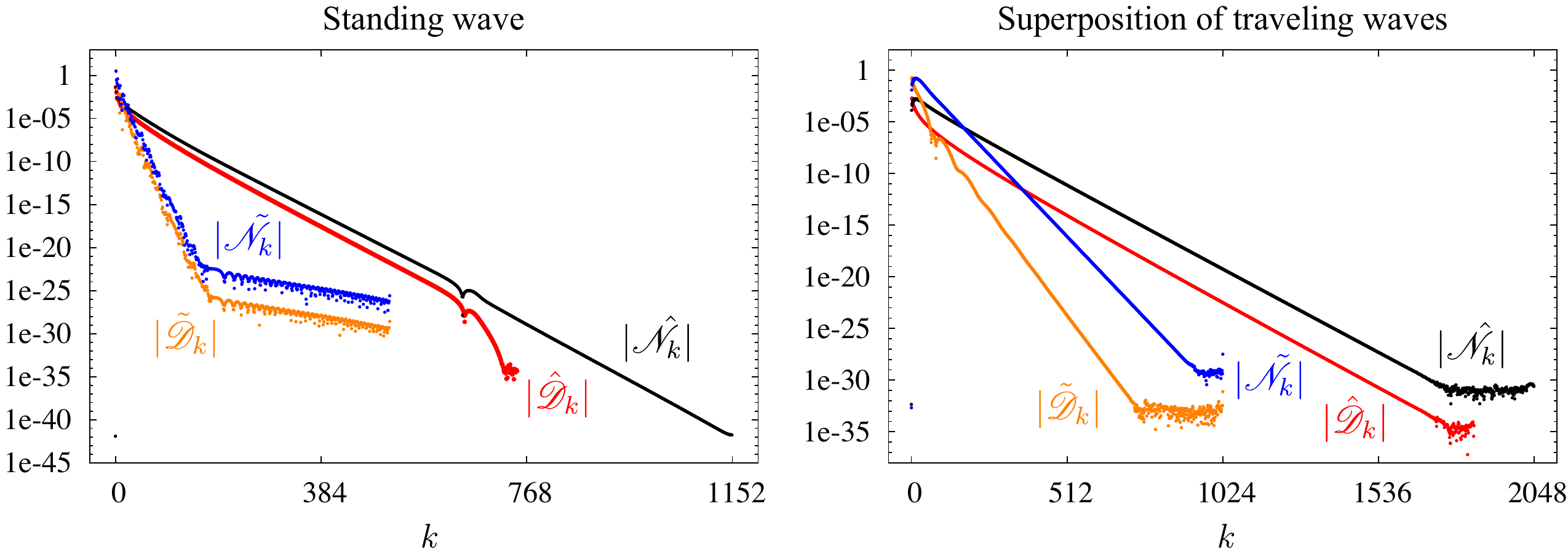}
\caption{\label{fig:afmDecay}
  Comparison of the coefficients in an expansion of the Dirichlet
  and Neumann data in an orthogonalized AFM basis versus a Fourier
  basis.
}
\end{figure}

The errors in the AFM and AFM$^*$ methods are nearly identical for
every choice of pseudo-inverse cutoff in the center panel of
Figure~\ref{fig:afmStokes}. Thus, the methods are equally effective at
computing the DNO for this example.  Beyond $k_\text{cutoff}=970$, the
error in the center panel reaches a plateau of $3\times10^{-28}$.
This is due to errors in the ``exact'' solution.  Since this plateau
is reached already at $k_\text{cutoff}=970$, it would appear that $K$
and possibly $M$ can be reduced while still achieving the target
accuracy of $10^{-28}$.  However, reducing $K$ leads to worse results
(not shown).  The reason is that the high-index columns of $A$ in
\begin{equation*}
  A_{j0} = \frac{1}{M}, \qquad
  \begin{pmatrix} A_{j,2k-1} \\ A_{j,2k} \end{pmatrix} =
  \frac{\sqrt{2}}{M}\,\frac{\cosh[k\eta(x_j)+h]}{\cosh[k\eta_\text{max}+h]}
  \begin{pmatrix} \cos kx \\[2pt] \sin kx \end{pmatrix}, \qquad
  (1\le k<K/2)
\end{equation*}
%
comprise a small but important part of the low-index singular vectors
$U$ of the singular value decomposition, $A=USV^T$.  A better way to
understand how large $K$ needs to be in order to achieve a prescribed
accuracy is to perform a QR factorization, $A=QR$, which is a
numerically robust way to perform a Gram-Schmidt orthogonalization of
the columns of $A$. We then define the ``discrete AFM transform''
\begin{equation*}
  \tilde\sN_0 = c_0, \qquad
  \tilde\sN_k = c_{2k-1}+ic_{2k}, \qquad
  \{c\}_{k=0}^{K-2} = \frac{1}{\sqrt{M}}Q^T\{\sN(x_j)\}_{j=0}^{M-1}
\end{equation*}
and, in Figure~\ref{fig:afmDecay}, compare it to the discrete Fourier
transform
\begin{equation*}
  \hat\sN_k = \frac{1}{M}\sum_{j=0}^{M-1} \sN(x_j)e^{-2\pi ijk}, \qquad
  0\le k\le M/2.
\end{equation*}
The red curves give the Fourier modes of the input Dirichlet data,
which are taken to specify the problem with infinite precision even
though they only describe the standing wave and traveling waves to
around 25 and 30 digits, respectively.  The black curves give the
Fourier modes of the boundary integral solutions, which were computed
with 212 and 106 bits of precision, respectively.  The orange and blue
curves give the ``AFM transform'' of the Dirichlet and Neumann
data. The observation that $K$ cannot be reduced significantly below
2048 for the traveling wave problem is seen clearly in the right
panel, which shows that $\tilde N_k$ does not reach roundoff error
until just below $k=1024$.

The fact that the AFM coefficients $\tilde D_k$ and $\tilde N_k$ decay
faster than the Fourier coefficients $\hat D_k$ and $\hat N_k$
suggests that the Gram-Schmidt version of the AFM basis is more
efficient at representing the Dirichlet and Neumann data of many
problems of physical interest (such as standing and traveling waves)
than the Fourier basis.  The abrupt change in slope of the orange and
blue curves in the left panel occurs when the AFM method begins to
resolve errors in the given Dirichlet data, which is taken to be exact
even though it only agrees with physical standing waves to 25 digits.
While it is interesting that the Gram-Schmidt version of the AFM basis
is more efficient than a Fourier basis for these problems, computing
the orthogonal AFM basis is expensive since it requires a great deal
of additional precision to avoid losing all significant digits during
the QR factorization process.

\section*{Conclusions}

We have shown that the Craig-Sulem expansion of the DNO operator, the
implicit formulation due to Ablowitz, Fokas and Musslimani, and its
dual, due to Ablowitz and Haut, can all be used to compute spectrally
accurate solutions of the Dirichlet-Neumann problem.  All three
methods involve ill-conditioned intermediate calculations when the
vertical crest-to-trough distance becomes large relative to the
inverse of the highest-frequency wave numbers involved in a Fourier
description of the wave profile $\eta(x)$ and the Dirichlet data
$\sD(x)$. However, this ill-conditioning can be tamed using extended
precision arithmetic.  Most importantly, $\eta(x)$ and $\sD(x)$ need
only be specified in double-precision to obtain double-precision
results for $\sN(x)$; extended precision is only required in
intermediate calculations.  While it is undesirable to work in
extended precision arithmetic, a spectrally accurate method requiring
64 digits of accuracy in intermediate calculations may require less
work than a more traditional 16-digit calculation using a 2nd or 4th
order method if high accuracy is desired.

Nevertheless, our goal in writing the paper was not to advocate the
use of these methods, but to explore their limits of applicability.
We originally thought that the AFM and AFM$^*$ would break down if
the potential $\phi$ cannot be extended analytically to a strip
containing the peak of the wave profile.  This was the motivation for
studying a potential $\phi$ with poles on the real axis and a wave
profile that dips below the poles but also extends above them.  It
is surprising that the AFM basis can still be used to approximate
$\sD$ in this case to arbitrary accuracy, and that the term-by-term
Neumann data turns out to be a good approximation of the correct
solution.  We also did not expect the CS expansion to converge for
realistic large-amplitude standing waves, but it was able to achieve
25 digits of accuracy at 150th order.  For the same problem, the AFM
and AFM$^*$ methods achieved 25 digits of accuracy in spite of
singular values dropping below $10^{-100}$.  For waves in shallow
water, the condition numbers are much more reasonable, making the
AFM methods more appealing.

As noted in \cite{NichollsReitich}, the CS expansion suffers from
significant cancellations requirements. Our numerical investigations 
confirm this behavior. We note that a suitable rearrangement of the 
series permits one to account for some of these cancellations,
though not all. One of the goals of the present work was to understand
whether this behavior is present in the AFM formulation and to what 
extent as the AFM method may be interpreted as a certain summation of
the CS expansion \cite{AblowitzHaut}. Our investigations reveal the AFM
method is similarly ill-conditioned. Indeed the 
rapid decay of singular values in the AFM (and AFM*) method mirrors the
cancellation properties observed in the CS expansion. Given the delicate 
nature of CS, AFM and AFM* we recommend the BIM, especially for 
two-dimensional flows. The TFE version of the CS expansion also performs
remarkably well without the need for extended precision arithmetic
in intermediate calculations.

Regarding extension to three dimensions, the CS and TFE expansions
work with little change.  Since the latter can be run effectively in
double-precision, it is a viable method.  Although the bulk fluid must
be discretized, the system of equations that must be solved decouples
into many one-dimensional boundary value problems rather than a large
system of equations such as arise in finite element methods.  The AFM
and AFM$^*$ methods will lead to large, highly ill-conditioned linear
systems that are easy to set up but nearly impossible to solve.
Performing the SVD or QR factorization of such a large matrix in
extended precision arithmetic would be extremely costly, and iterative
methods such as GMRES will not converge when the condition number
grows to $10^{15}$ or higher.  By contrast, the boundary integral
method yields condition numbers close to 1; hence, GMRES converges in
just a few iterations.  Although it is difficult to implement boundary
integral methods in three dimensions due to the slowly decaying
lattice sums involved, techniques such as Ewald summation
\cite{krasny:ewald,annaKarin,siegel} are available to accelerate
convergence. 

\section*{Acknowledgments}
JW was supported in part by the Director, Office of Science,
Computational and Technology Research, U.S.  Department of Energy
under Contract No.~DE-AC02-05CH11231, and by the National Science
Foundation through grant DMS-0955078.  Any opinions, findings, and
conclusions or recommendations expressed in this material are those of
the authors and do not necessarily reflect the views of the funding
sources.

\bibliographystyle{plainnat}

\begin{thebibliography}{33}
\providecommand{\natexlab}[1]{#1}
\providecommand{\url}[1]{\texttt{#1}}
\expandafter\ifx\csname urlstyle\endcsname\relax
  \providecommand{\doi}[1]{doi: #1}\else
  \providecommand{\doi}{doi: \begingroup \urlstyle{rm}\Url}\fi

\bibitem[Ablowitz and Haut(2008)]{AblowitzHaut}
M.~J. Ablowitz and T.~S. Haut.
\newblock {Spectral formulation of the two fluid Euler equations with a free
  interface and long wave reductions}.
\newblock \emph{Analysis and Applications.}, 6:\penalty0 323--348, 2008.

\bibitem[Ablowitz et~al.(2006)Ablowitz, Fokas, and Musslimani]{AFM}
M.~J. Ablowitz, A.~S. Fokas, and Z.~H. Musslimani.
\newblock {On a new non-local formulation of water waves}.
\newblock \emph{J. Fluid Mech.}, 562:\penalty0 313--343, 2006.

\bibitem[Ambrose et~al.(2013)Ambrose, Siegel, and Tlupova]{siegel}
D.~M. Ambrose, M.~Siegel, and S.~Tlupova.
\newblock A small-scale decomposition for {3D} boundary integral computations
  with surface tension.
\newblock \emph{J. Comput. Phys.}, 247:\penalty0 168--191, 2013.

\bibitem[Baker and Nachbin(1998)]{baker:nachbin:98}
G.~Baker and A.~Nachbin.
\newblock Stable methods for vortex sheet motion in the presence of surface
  tension.
\newblock \emph{SIAM J. Sci. Comput.}, 19\penalty0 (5):\penalty0 1737--1766,
  1998.

\bibitem[Baker and Xie(2011)]{baker10}
G.~R. Baker and C.~Xie.
\newblock Singularities in the complex physical plane for deep water waves.
\newblock \emph{J. Fluid Mech.}, 685:\penalty0 83--116, 2011.

\bibitem[Baker et~al.(1982)Baker, Meiron, and Orszag]{baker:82}
G.~R. Baker, D.~I. Meiron, and S.~A. Orszag.
\newblock Generalized vortex methods for free-surface flow problems.
\newblock \emph{J. Fluid Mech.}, 123:\penalty0 477--501, 1982.

\bibitem[Bruno and Reitich(1992)]{brunoReitich}
O.P. Bruno and F.~Reitich.
\newblock Solution of a boundary value problem for the {Helmholtz} equation via
  variation of the boundary into the complex domain.
\newblock \emph{Proc. Royal Soc. Edinburgh}, 122A:\penalty0 317--340, 1992.

\bibitem[Canuto et~al.(1988)Canuto, Hussaini, Quarteroni, and Zang]{canuto:88}
C.~Canuto, M.~Y. Hussaini, A.~Quarteroni, and T.~A. Zang.
\newblock Spectral methods in fluid dynamics.
\newblock 1988.

\bibitem[Cohen and Kundu(2004)]{Kundu}
I.~M. Cohen and P.~K. Kundu.
\newblock \emph{{Fluid Mechanics}}.
\newblock Academic Press, 2004.

\bibitem[Craig and Sulem(1993)]{CraigSulem}
W.~Craig and C.~Sulem.
\newblock Numerical simulation of gravity waves.
\newblock \emph{J. Comp. Phys.}, 108:\penalty0 73--83, 1993.

\bibitem[Craig et~al.(2005)Craig, Guyenne, Nicholls, and Sulem]{Craigetal}
W.~Craig, P.~Guyenne, D.~P. Nicholls, and C.~Sulem.
\newblock {Hamiltonian long-wave expansions for water waves over a rough
  bottom}.
\newblock \emph{Proc. R. Soc. A}, 461:\penalty0 839--873, 2005.

\bibitem[Deconinck and Oliveras(2011)]{OliverasDeconinck}
B.~Deconinck and K.~Oliveras.
\newblock {The instability of periodic surface gravity waves}.
\newblock \emph{J. Fluid Mech.}, 675:\penalty0 141--167, 2011.

\bibitem[Duan and Krasny(2000)]{krasny:ewald}
Z.-H. Duan and R.~Krasny.
\newblock An {Ewald} summation based multipole method.
\newblock \emph{J. Chem. Phys.}, 113\penalty0 (9):\penalty0 3492--3495, 2000.

\bibitem[Dyachenko et~al.(1996)Dyachenko, Kuznetsov, Spector, and
  Zakharov]{dyachenko1996}
A.I. Dyachenko, E.A. Kuznetsov, M.D. Spector, and V.E. Zakharov.
\newblock Analytical description of the free surface dynamics of an ideal fluid
  (canonical formalism and conformal mapping).
\newblock \emph{Phys. Letters A}, 221:\penalty0 73 -- 79, 1996.

\bibitem[Fousse et~al.(2007)Fousse, Hanrot, Lef\`evre, P\'elissier, and
  Zimmermann]{mpfr}
Laurent Fousse, Guillaume Hanrot, Vincent Lef\`evre, Patrick P\'elissier, and
  Paul Zimmermann.
\newblock {MPFR}: A multiple-precision binary floating-point library with
  correct rounding.
\newblock \emph{{ACM} Transactions on Mathematical Software}, 33\penalty0
  (2):\penalty0 13:1--13:15, June 2007.
\newblock URL \url{http://doi.acm.org/10.1145/1236463.1236468}.

\bibitem[Krasny(1986)]{krasny:86}
R.~Krasny.
\newblock Desingularization of periodic vortex sheet roll-up.
\newblock \emph{J. Comput. Phys.}, 65:\penalty0 292--313, 1986.

\bibitem[Lannes(2005)]{Lannes}
D.~Lannes.
\newblock Well-posedness of the water-wave equations.
\newblock \emph{J. Amer. Math. Soc.}, 18:\penalty0 605--654, 2005.

\bibitem[Lindbo and Tornberg(2012)]{annaKarin}
D.~Lindbo and A.-K. Tornberg.
\newblock Fast and spectrally accurate {Ewald} summation for 2-periodic
  electrostatic systems.
\newblock \emph{J. Chem. Phys.}, 136:\penalty0 164111, 2012.

\bibitem[Longuet-Higgins and Cokelet(1976)]{lh76}
M.~S. Longuet-Higgins and E.~D. Cokelet.
\newblock The deformation of steep surface waves on water. {I}. a numerical
  method of computation.
\newblock \emph{Proc. Royal Soc. A}, 350:\penalty0 1--26, 1976.

\bibitem[Mercer and Roberts(1992)]{mercer:92}
G~N Mercer and A~J Roberts.
\newblock {Standing waves in deep water: Their stability and extreme form}.
\newblock \emph{Phys. Fluids A}, 4\penalty0 (2):\penalty0 259--269, 1992.

\bibitem[Mercer and Roberts(1994)]{mercer:94}
G~N Mercer and A~J Roberts.
\newblock The form of standing waves on finite depth water.
\newblock \emph{Wave Motion}, 19:\penalty0 233--244, 1994.

\bibitem[Muskhelishvili(1992)]{muskhelishvili}
N.~I. Muskhelishvili.
\newblock \emph{Singular Integral Equations, 2nd Edition}.
\newblock Dover, New York, 1992.

\bibitem[Nicholls and Reitich(2001)]{NichollsReitich}
D.~P. Nicholls and F.~Reitich.
\newblock A new approach to analyticity of {D}irichlet-{N}eumann operators.
\newblock \emph{Proc. Roy. Soc. Edin.: Sec. A Mathematics}, 131:\penalty0
  1411--1433, 2001.

\bibitem[Nicholls and Reitich(2006)]{NichollsReitich06}
D.~P. Nicholls and F.~Reitich.
\newblock Stable, high-order computation of traveling water waves in three
  dimensions.
\newblock \emph{European J. Mech. B/Fluids}, 25:\penalty0 406--424, 2006.

\bibitem[Oliveras et~al.(2012)Oliveras, Vasan, Deconinck, and Henderson]{OVDH}
K.~Oliveras, V.~Vasan, B.~Deconinck, and D.~Henderson.
\newblock Recovering the water-wave profile from pressure measurements.
\newblock \emph{{SIAM J. Appl. Math.}}, 72:\penalty0 897--918, 2012.

\bibitem[Rycroft and Wilkening(2013)]{rycroft:13}
C.~H. Rycroft and J.~Wilkening.
\newblock Computation of three-dimensional standing water waves.
\newblock \emph{J. Comput. Phys.}, 255:\penalty0 612--638, 2013.

\bibitem[Smith and Roberts(1999)]{smith:roberts:99}
D.~H. Smith and A.~J. Roberts.
\newblock Branching behavior of standing waves --- the signatures of resonance.
\newblock \emph{Phys. Fluids}, 11:\penalty0 1051--1064, 1999.

\bibitem[Vasan and Deconinck(2013)]{vasanDeconinck}
V.~Vasan and B.~Deconinck.
\newblock The inverse water wave problem of bathymetry detection.
\newblock \emph{Journal of Fluid Mechanics}, 714:\penalty0 562--590, 2013.

\bibitem[Wilkening(2014)]{trav:stand}
J.~Wilkening.
\newblock Traveling-standing water waves.
\newblock 2014.
\newblock (submitted).

\bibitem[Wilkening et~al.(2014)Wilkening, Cerfon, and Landreman]{plasma2}
J.~Wilkening, A.~Cerfon, and M.~Landreman.
\newblock Projected dynamics of kinetic equations with energy diffusion in
  spaces of orthogonal polynomials.
\newblock 2014.
\newblock (submitted).

\bibitem[Wilkening(2011)]{water1}
Jon Wilkening.
\newblock Breakdown of self-similarity at the crests of large-amplitude
  standing water waves.
\newblock \emph{Phys. Rev. Lett.}, 107:\penalty0 184501, Oct 2011.
\newblock \doi{10.1103/PhysRevLett.107.184501}.

\bibitem[Wilkening and Yu(2012)]{water2}
Jon Wilkening and Jia Yu.
\newblock Overdetermined shooting methods for computing standing water waves
  with spectral accuracy.
\newblock \emph{Comput. Sci. Disc.}, 5\penalty0 (1):\penalty0 014017, 2012.
\newblock \doi{10.1088/1749-4699/5/1/014017}.

\bibitem[Zakharov(1968)]{Zakharov}
V.~E. Zakharov.
\newblock {Stability of periodic waves of finite amplitude on the surface of a
  deep fluid}.
\newblock \emph{Zhurnal Prikladnoi Mekhaniki i Tekhnicheskoi Fiziki},
  8:\penalty0 86--94, 1968.

\end{thebibliography}

\end{document}